\documentclass{article}
\usepackage{amssymb}
\usepackage{amsfonts}
\usepackage{amsmath,amsthm,amssymb,latexsym,amscd}
\usepackage{epsf}
\usepackage[pdftex]{graphicx}
\usepackage{graphics}

\newtheorem{lemma}{Lemma}
\newtheorem{theorem}{Theorem}
\begin{document}

\title{{\Large \textbf{Optimal Local Approximation Spaces for Generalized
Finite Element Methods with Application to Multiscale Problems\thanks{%
This work is supported by grants: NSF DMS-0807265 and AFOSR
FA9550-08-1-0095.}}}}
\author{Ivo Babuska $^{\mbox{\tiny 1}}$ \and Robert Lipton $^{\mbox{\tiny 2}%
} $}
\date{}
\maketitle

\begin{abstract}
The paper addresses a numerical method for solving second order elliptic
partial differential equations that describe fields inside heterogeneous
media. The scope is general and treats the case of rough coefficients, i.e.
coefficients with values in $L^\infty(\Omega)$. This class of coefficients
includes as examples media with micro-structure as well as media with
multiple non-separated length scales. The approach taken here is based on
the the generalized finite element method (GFEM) introduced in \cite{107},
and elaborated in \cite{102}, \cite{103} and \cite{104}. The GFEM is
constructed by partitioning the computational domain $\Omega$ into a
collection of preselected subsets $\omega _{i},i=1,2,..m$ and constructing
finite dimensional approximation spaces $\Psi _{i}$ over each subset using
local information. The notion of the Kolmogorov $n$-width is used to
identify the optimal local approximation spaces. These spaces deliver local
approximations with errors that decay almost exponentially with the degrees
of freedom $N_{i}$ in the energy norm over $\omega_i$. The local spaces $%
\Psi _{i}$ are used within the GFEM scheme to produce a finite dimensional
subspace $S^N$ of $H^{1}(\Omega )$ which is then employed in the Galerkin
method. It is shown that the error in the Galerkin approximation decays in
the energy norm almost exponentially ( i.e., super-algebraicly) with respect
to the degrees of freedom $N$. When length scales ``separate'' and the
microstructure is sufficiently fine with respect to the length scale of the
domain $\omega_i$ it is shown that homogenization theory can be used to
construct local approximation spaces with exponentially decreasing error in
the pre-asymtotic regime.
\end{abstract}

\baselineskip=0.9\normalbaselineskip
\vspace{-3pt}

\begin{center}
{\footnotesize $^{\mbox{\tiny\rm 1}}$ICES and Department of Aerospace
Engineering, University of Texas, Austin, TX \\[0pt]
78712, USA. email: babuska\symbol{'100}ices.utexas.edu\\[3pt]
$^{\mbox{\tiny\rm 2}}$Department of Mathematics, Louisiana State University,
Baton Rouge, LA\\[0pt]
70803, USA. email: lipton\symbol{'100}math.lsu.edu\\[3pt]
}
\end{center}

\setcounter{equation}{0} \setcounter{theorem}{0} \setcounter{lemma}{0}

\section{Introduction}

Large multi-scale systems such as airplane wings and wind turbine blades are
built from fiber reinforced composites and exhibit a cascade of substructure
spread across several length scales. These and other large composite
structures are now seeing extensive use in transportation, energy, and
infrastructure. The importance of accurate numerical simulation is ever
increasing due to the high cost of experimental testing of large structures
made from heterogeneous materials. The computational modeling of such
heterogeneous structures is a very large problem that requires the use of
parallel computers. 
In order for a numerical method to be adequate it must
be able to utilize many local computations  performed
independently on single processors or clusters of processors of reasonable
size. 
The approach taken here is based on the Generalized Finite Element
Method (GFEM) introduced in \cite{107}, and elaborated in \cite{102}, \cite%
{103} , \cite{104}. It is a partition of unity method \cite{107} which
utilizes the results of many independent and local computations carried out
across the computational domain. The GFEM is constructed by partitioning the
computational domain $\Omega $ into to a collection of preselected subsets $%
\omega _{i},i=1,2,..m$ and constructing finite dimensional approximation
spaces $\Psi _{i}$ over each subset using local information. The specific
way in which the partition is carried out is special to this method \cite{107}, 
\cite{102} and the details are discussed in
section \ref{GFEM}. Since each space $\Psi _{i}$ is computed independently
the full \textquotedblleft global\textquotedblright\ solution is obtained by
solving a global (macro) system which is an order of magnitude smaller than
the system corresponding to a direct application the finite element method
to the full structure. The GFEM approach provides an opportunity for the
significant reduction of the computational work involved in the numerical
modeling of large heterogeneous problems.

In this article we show how to achieve optimal accuracy within the GFEM
approach. The key point is to note that the approximation error of the GFEM
is controlled by the corresponding approximation error of the local
approximation spaces $\Psi_i,\,i=1,\ldots,m$, see section \ref{GFEM}.
Therefore the goal is to identify optimal local approximation spaces. 
Our approach is naturally guided by the notion of the Kolmogorov $n$-width 
\cite{Pincus} which measures the ability of an increasing sequence of finite
dimensional subspaces of a prescribed Banach space $B$ to approximate any
element inside $B$, see section \ref{sectexponentoverview}. Using the
solution of the spectral problem associated with the $n$-width we are able
to identify a new class of approximation spaces $\Psi _{i}$. We show that
these finite dimensional spaces are able to approximate the solution on $%
\omega _{i}$ with errors that decay almost exponentially with the degrees of
freedom $N_{i}$ in the energy norm. The overall method for constructing the
local approximations is general and applies to subdomains $\omega_i$
belonging to the interior of the computational domain as well as those that
intersect the boundary of the computational domain. The optimal local
approximation spaces are identified for interior subdomains in section 
\ref{locinterior} and for those touching the boundary of the computational
domain in section \ref{locbdry}.

The optimal local spaces $\Psi _{i}$ are then combined within the GFEM
scheme to produce a finite dimensional GFEM subspace $S^N$ of $H^{1}(\Omega
) $ which is employed in the Galerkin method. 
It is shown that the
corresponding error in the Galerkin approximation decays in the energy norm
almost exponentially ( i.e., super-algebraicly) with respect to the degrees
of freedom $N$ see section \ref{Macroexponential}. 
In section \ref%
{implementation} we discuss the main issues involved in the implementation
as well as estimates of the computational work associated with this
numerical approach.

We show how to construct nearly optimal local approximation spaces using
homogenized coefficients when the subdomain $\omega_i$ is sufficiently large
with respect to the length scale of the heterogeneity. The homogenization
limit for the $n$-width and the optimal approximation space is identified in
section \ref{nwidthhomog}. This identification is established within the
general homogenization context described by $H$-convergence and $G$%
-convergence, \cite{MuratTartar}, \cite{7}. These results are applied to
heterogeneous media with micro-structure that has \emph{uniformly fine}
variation with respect to the length scale of the domains $\omega_i$. Here a 
\emph{uniformly fine} microstructure is defined to be one that can be
identified as belonging to a sequence of microstructures characterized by a
sequence of length scales $\epsilon=\frac{1}{k},\,k=1,2,\ldots$ and
coefficients $\{A^\epsilon\}_{\epsilon>0}$ that converge to a homogenization
limit described by a matrix of constant coefficients $A$. For this case
we provide examples that illustrate how to construct local approximation
spaces with errors that decay exponentially in the pre-asymptotic regime,
see section \ref{homogpre}. The examples corroborate the exponentially
decreasing error observed in the numerical simulations for finely mixed
dispersed inclusions carried out in \cite{108}.

The homogenization theory developed here provides motivation for some rules
of thumb for choosing the size of subdomains $\omega_i$ in the
implementation. Here the size is chosen large relative to the local length
scale of the heterogeneity but small enough such that the heterogeneity is
statistically uniform within it. The specific details are presented in
section \ref{homogpre}.

We conclude noting that there is now a large and rapidly growing literature
devoted to the numerical analysis of multi-scale media. Several contemporary
mathematically based approaches include the Multiscale FEM \cite{Hou1}, \cite%
{Hou2}, global changes of coordinates for upscaling porous media flows \cite%
{Efendiev1}, \cite{Efendiev2} , the heterogeneous multiscale methods (HMM) 
\cite{WE1},\cite{WE2}, \cite{Engq}, an adaptive coarse scale - fine scale
projection method \cite{Nolen}, numerical homogenization methods for $%
L^\infty$ coefficients based on harmonic coordinates and elliptic
inequalities \cite{OhwadiZhang1}, \cite{OhwadiZhang2}, \cite{OhwadiBerlyand}%
, subgrid upscaling methods \cite{Arbogast1}  and global
Galerkin projection schemes for problems with $L^\infty$ coefficients and
homogeneous Dirichlet boundary data \cite{melink}.

\subsection{Problem formulation}

\label{problemformulation} Let $\ \Omega \in \mathbb{R}^{d}$ be a bounded domain
with $C^{1}$ smooth boundary $\partial \Omega $. In this article we consider
the elliptic differential equation 
\begin{equation}
-div(A(x)\nabla u(x))=f(x),\hbox{    $\forall x\in \Omega $}
\label{1.1}
\end{equation}%
with either Neumann boundary conditions prescribed on the boundary $\partial
\Omega $ given by 
\begin{equation}
n\cdot A\nabla u(x)=g,\hbox{   $x\in
\partial \Omega $}
\label{1.2a}
\end{equation}%
where $n$ is the outer unit normal vector or Dirichlet boundary conditions 
\begin{equation}
u(x)=q,\hbox{   $x\in \partial
\Omega$}.
\label{Dirichletbc}
\end{equation}%
\noindent In forthcoming work we will address the case of non-smooth boundaries and the case where both
Dirichlet and Neumann boundary conditions are prescribed on different parts
of the boundary.

We assume that $A(x)$ is a $d\times d$ symmetric matrix with measurable
coefficients $a_{i,j}(x)$ $\in L^{\infty }(\Omega )$ satisfying the standard
coercivity condition 
\begin{eqnarray}  \label{1.3a}
\alpha ^{\ast }(x)\mid \mathbf{v}\mid ^{2}\leq \mathbf{v}^{t}A(x)\mathbf{%
v\leq \beta }^{\ast }(x)\mathbf{\mid v\mid }^{2}\mathbf{,}\forall x\in \Omega
\end{eqnarray}
where $\mathbf{v\in R}^{d}$ is an arbitrary vector and 
\begin{eqnarray}  \label{1.3b}
0<\alpha \leq \alpha ^{\ast }(x)\leq \beta ^{\ast }(x)\leq \beta <\infty , %
\hbox{  $\forall x\in \Omega.$}
\end{eqnarray}
For future reference we will denote this class of coefficients by $\mathfrak{%
C}$. Here we will suppose that $f\in H^{k}(\Omega ),k\geq 0,$ $g\in
H^{-1/2}(\partial \Omega ),$ together with the consistency condition $%
\int_{\partial \Omega }gds+\int_{\Omega }fdx=0$\ for the Neumann boundary
condition \eqref{1.2a} and that the Dirichlet data $q$ is an element of $%
H^{1/2}((\partial \Omega )$.
The weak solution belonging to $u\in H^{1}(\Omega )$ exists and for Dirichlet boundary
conditions it is unique and for Neumann boundary conditions is unique up to
additive constant. 

%Because we have assumed only that the coefficients are
%measurable we can find  rough coefficients 
%for which solutions $u\notin H^{1+\varepsilon }(\Omega )$ for any $%
%\epsilon >0$.

In this work we investigate the accuracy of the GFEM for the approximation
of the exact solution $u_{0}$. Here the objective is
to find an approximate solution $u_{\tau }\in H^{1}(\Omega )$ for which 
\begin{equation*}
\parallel u_{0}-u_{\tau }\parallel _{\mathcal{E}(\Omega )}\leq \tau
\end{equation*}%
where 
\begin{equation}
\parallel u\parallel _{\mathcal{E}(\Omega )}=\left( \int_{\Omega }A\nabla
u\cdot \nabla u\,dx\right) ^{1/2}
\label{energynorm}
\end{equation}%
is the energy norm. In this treatment we address the scalar problem noting
that the ideas used in the approach presented here apply with out any
modification to second order elliptic systems including the system of
linear elasticity.

\subsection{A typical example: fiber reinforced composites}

\begin{figure}[th]
%[t]
\centering
\scalebox{0.7}{\includegraphics{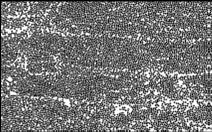}}
\caption{Fiber-reinforced composite.}
\label{Plate}
\end{figure}

\begin{figure}[th]
%[t]
\centering
\scalebox{0.4}{\includegraphics{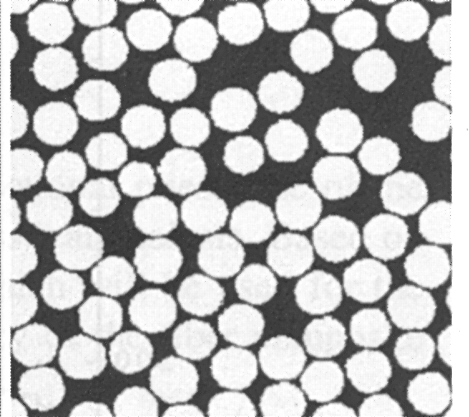}}
\caption{Microstructure.}
\label{small_sample}
\end{figure}

For the purposes of this article we have chosen to work with coefficient
matrices belonging to $L^{\infty }(\Omega )$ subject to standard coercivity
and boundedness conditions. This choice reflects our intention to describe
generic situations for which there can be several non-separated length
scales of variation inside the heterogeneous media. In this treatment we
shall assume that the coefficient matrix describing the media is known. To
fix ideas we discuss the problem of determining the coefficient matrix $A$
associated with a sample of fiber reinforced composite material described in 
\cite{105}. A fiber reinforced material consists of two components the fiber
and the host material commonly referred to as the matrix. The material shown
here is taken from the center of a composite plate of 36 plies, divided into
9 groups each containing 4 plies. The orientation of the fibers alternates
between 0$^{o}$ and 90$^{o}$ from group to group. The plys are HTA/8376
unidirectional prepreg fiber composites produced by Ciba-Geigy. The sample
is a rectangular plate of length 300mm and width 140 mm. The nominal ply
thickness is 130 $\mu m$. Figure \ref{Plate} shows the cross section of 
a group of four plies consisting of 16275 fibers. The cross section
of each fiber is roughly circular with a fiber diameter of about 7$\mu m$.
It is clear from Figure \ref{Plate} that the material coefficients have
variation across several length scales. These include \textquotedblleft
matrix rich\textquotedblright\ zones between the plys as well as variation
in fiber alignment between groups of plys. At the smallest length scale
Figure \ref{small_sample} shows that the coefficients are piece wise
constant taking one value in the fiber and a different value in the matrix.
This composite sample has been mapped in the study \cite{105} and 
the relative fiber positions are known
exactly and so the coefficient matrix can be determined. 
In general it is impossible to record the relative location of every fiber
for an entire structure made from a composite material. Hence some
stochastic information needs to be extracted from the structure and used to
describe the coefficient matrix. An investigation of the different
stochastic data that can be obtained from composite samples is taken up in \cite{105}. As mentioned
earlier we will assume that the problem is deterministic with well defined
coefficients. Future work will address the problem of determination of the
approximation error for stochastically defined coefficients.

%%%%%%%%%%%
%Finally in these problems the functions $f,g,q$ in practice are very, smooth
%analytic or piecewise analytic so that the exponential convergence rate can
%be obtained.

\setcounter{equation}{0} \setcounter{theorem}{0} \setcounter{lemma}{0}

\section{The generalized finite element method (GFEM)}

\label{GFEM}

The GFEM is based on the partition of unity method (PUM) and is introduced
in \cite{107} as a method for the numerical solution of elliptic PDE with
rough coefficients. The method is further elaborated and extended to other
application areas in the works of \cite{102}, \cite{103}, \cite{101} and 
\cite{104}. We briefly summarize the main ideas and results of GFEM. For
more details see \cite{102}. We recall the scalar Dirichlet and Neumann
boundary value problems defined in section \ref{problemformulation}. For the
Dirichlet problem define the hyperplane $H_{q}^{1}(\Omega )=\{u\in
H^{1}(\Omega )\mid u=q(x)$ on $\partial \Omega \}$ and set $%
H_0^1(\Omega)=H_q^1(\Omega)$ for $q=0$. The weak solution of the Dirichlet
problem $u_{0}\in H_{q}^{1}(\Omega )$ satisfies 
\begin{eqnarray}  \label{2.1}
B(u_{0},v)=F(v), 
\end{eqnarray}
for all $v\in
H_{0}^{1}(\Omega )$ where 
\begin{eqnarray}  \label{2.2}
B(u,v)=\int_{\Omega } A(x)\nabla u \cdot \nabla v\,dx,&\hbox{  and  }&
F(v)=\int_{\Omega }f\,v dx.
\end{eqnarray}

For the Neumann boundary value problem the solution $u_{0}\in H^{1}(\Omega )$
satisfies \eqref{2.1} for all $v\in H^{1}(\Omega )$ and 
\begin{equation*}
F(v)=\int_{\Omega }fvdx+\int_{\partial \Omega }gvds
\end{equation*}%
For future reference we note that the energy norm \eqref{energynorm} is
given by $\parallel u\parallel _{\mathcal{E}(\Omega )}=(B(u,u))^{1/2}$.

We now recall that for a $N$ dimensional space $%
S(\Omega )\subset H^{1}(\Omega )$ that the associated
Galerkin solution of the Neumann problem \ $u^{S}\in S(\Omega )\subset
H^{1}(\Omega )$ satisfies $B(u^{S},v)=F(v)$ for all $v\in S(\Omega )$. For a fixed
tolerance $\varepsilon >0$ if there is a $\psi \in S(\Omega )$ such that $\parallel
u_{0}-\psi \parallel _{\scriptscriptstyle{\mathcal{E}(\Omega )}}\leq
\varepsilon $ then it is clear that the Galerkin solution satisfies $\parallel
u_{0}-u^{S}\parallel _{\scriptscriptstyle{\mathcal{E}(\Omega )}}\leq
\varepsilon $. 
A similar
conclusion holds for the Galerkin solution of the Dirichlet problem. 
Let $S_{0}(\Omega )$ be an $N$
dimensional subspace of $H_{0}^{1}(\Omega )$ and set $S_{q}(\Omega )=S_{0}(\Omega )\oplus \Phi _{0}$ where $\Phi
_{0}$ is a particular function belonging to $H_{q}^{1}$($\Omega )$. Then
the Galerkin solution of the Dirichlet problem $u^{S}\in S_{q}(\Omega )$
satisfies $B(u^{S},v)=F(v)$ for $\forall v\in S_{0}(\Omega ).$ Again
it is clear that if there exists a $\psi \in S_{q}(\Omega )$ such
that $\parallel u_{0}-\psi \parallel _{\scriptscriptstyle{\mathcal{E}(\Omega
)}}\leq \varepsilon $ then  $\parallel u_{0}-u^{S}\parallel
_{\scriptscriptstyle{\mathcal{E}(\Omega )}}\leq \varepsilon $. 
Hence the main task in constructing a Galerkin numerical
solution is the selection of $S(\Omega )$
for the Neumann problem and the selection of the hyperplane $S_{q}(\Omega )$ for the Dirichlet problem.

We introduce the Galerkin approximation delivered by GFEM and discuss the associated approximation error. 
Since most physical situations are described by  Neumann boundary conditions we address this case and note 
that the Dirichlet case follows identical lines.    Let $\{\mathcal{O}_{i}\}_{i=1}^m$, be a collection of 
open sets that cover the computational domain, i.e., $\Omega\subset\cup _{i=1}^{m}\mathcal{O} _{i}$, and 
we introduce the partition of unity subordinate to the open cover denoted by $\phi_{i}\in C^{1}(\mathcal{O}_{i})$, $i=1,...m$. 
We relabel interior sets $\mathcal{O}_i\subset\Omega$ 
as $\omega_i=\mathcal{O}_i$ and for sets $\mathcal{O}_i$ that intersect the boundary of $\Omega$ we write
$\omega_i=\mathcal{O}_i\cap\Omega$. Here $\cup_{i=1}^m\omega_i=\Omega$ and we assume that each point $x\in \Omega $
belongs to at most $\kappa$ subdomains $\omega _{i}$. The functions $\phi_{i}$, $i=1,\ldots,m$ have the following properties
\begin{eqnarray}
0 \leq \phi _{i}\leq 1,&&\hbox{  $i=1,\ldots,m,$}  \label{2.3a} \\
\phi _{i}(x)=0, &&%
\hbox{  for $x\in \Omega\backslash \omega _{i},$
$i=1,..m$} \label{2.3b}\\
\sum_{i=1}^{m}\phi _{i}(x)=1, &&\hbox{   $\forall x\in \Omega $},\label{2.3c} \\
\max_{x\in \Omega }\mid \phi _{i}(x)\mid &\leq &C_{1},\,\,\,\,i=1,...m \label{2.3d}\\
\max_{x\in \Omega }\mid \nabla \phi _{i}(x)\mid &\leq &\frac{C_{2}}{%
diam(\omega _{i})},\,\,\,\,i=1,...m,\label{2.3e}
\end{eqnarray}%
where the constants $C_{1}$ and $C_{2}$ are positive and bounded. Here $diam(\omega _{i})$ denotes the diameter
of $\omega _{i}.$ %Conditions, \eqref{2.3a}, \eqref{2.3b} and \eqref{2.3c} state that $%
%\{\phi _{i}\}$ is a partition of unity on $\Omega $. A typical example is
%given by a square mesh together with the element-wise bilinear basis
%functions (the standard hat functions) in two dimensions.

Next we introduce local approximation spaces associated with each $\omega
_{i}$. For this case let $\Psi _{i}$ be a finite
dimensional subspace of $H^{1}(\omega _{i})$ of dimension $N_{i}$.
The trial and test spaces for the GFEM are constructed from the local approximation spaces and are defined by

\begin{equation}
S(\Omega )=\left\{ \psi ;\,\psi =\sum_{i=1}^{m}\sum_{j=1}^{N_{i}}\phi
_{i}\xi _{i}^{[j]},\hbox{ where }\xi _{i}^{[j]}\in \Psi _{i},\,\right\}.
\label{GFEMBASIS}
\end{equation}
\noindent Note that although $\xi
_{i}^{[j]}$ only belongs to $H^{1}(\omega _{i})$ the PUM construction ensures  that $S(\Omega )\subset H^{1}(\Omega )$. Here the dimension $N$ of $%
S(\Omega )$ is given by $N=\sum_{i=1}^{m}N_{i}$.

An analogous approach using local approximation spaces is used for the Dirichlet problem. Here the only difference is when $\overline{%
\omega }_{j}\cap \partial \Omega \neq 0$. For this case if $q=0$ then $\Psi _{j}\subset
H^{1}(\omega _{j})\cap H_{0}^{1}(\Omega )$ and if $q\neq 0$ then $\Psi
_{j}\subset (H^{1}(\omega _{j})\cap H_{0}^{1}(\Omega ))\oplus \Phi _{j}$
where $\Phi _{j}$ is a particular function belonging to $H^{1}(\omega _{j})\cap
H_{q}^{1}(\Omega )$.

It is now clear from the formulation that the approximation error of the Galerkin numerical solution
for the GFEM is tied to the accuracy of the local approximation spaces.
With this in mind we state the following approximation theorem \cite{103} for the Neumann problem.

\begin{theorem}
\label{globallocal} Suppose that there exists $\zeta _{i}\in \Psi _{i}$, $%
i=1,\ldots ,m$ with 
\begin{eqnarray}
\parallel u_{0}-\zeta _{i}\parallel _{L^2_\ast(\omega _{i})} &=&(\int_{\omega
_{i}}\beta ^{\ast }(u_{0}-\zeta _{i})^{2}dx)^{1/2}\leq \varepsilon _{1}(i)
\label{2.4b} \\
\parallel u_{0}-\zeta _{i}\parallel _{\mathcal{E}(\omega _{i})} &\leq
&\varepsilon _{2}(i)
\end{eqnarray}%
and let 
\begin{equation*}
\zeta (x)=\sum_{i=1}^{m}\zeta _{i}(x)\phi _{i}(x)
\end{equation*}%
then $\zeta (x)\in H^{1}(\Omega )$ and 
\begin{eqnarray}
\parallel u_{0}-\zeta \parallel _{L^{2}_\ast(\Omega )} &\leq &C_{1}(\sum_{i=1}^{m}(\varepsilon _{1}(i))^{2})^{1/2}  \label{2.6c} \\
\parallel u_{0}-\zeta \parallel _{\mathcal{E}(\Omega )} &\leq &(C_{2}^{2}\sum_{i=1}^{m}(\varepsilon _{1}^{2}(i)/diam^{2}(\omega
_{i}))+C_{1}^{2}\sum_{i=1}^{m}\varepsilon _{2}^{2}(i))^{1/2},\label{2.6d}
\end{eqnarray}%
where $C_{1}$ and $C_{2}$\ are as in \eqref{2.3d} and \eqref{2.3e}
respectively.

Moreover if all local spaces $\Psi _{i}$ contain the subspace of constant functions
then, for an appropriate choice of constant, the first term on the right hand side of \eqref{2.6d} can be omitted as
it is majorized by the second term. 
\end{theorem}

We remark that Theorem \ref{globallocal} also holds true for Dirichlet boundary conditions.
For this case the only modifications are the ones previously discussed for the subdomains $\omega_j$ that intersect the boundary. Theorem \ref{globallocal}  shows that the proper selection of the local
approximation spaces $\Psi _{i}$ are essential for obtaining optimal
accuracy. In the next section we identify local spaces $\Psi _{i}$ that
deliver a nearly exponential rate of convergence with respect to the
dimension of the trial space.

\setcounter{equation}{0} \setcounter{theorem}{0} \setcounter{lemma}{0}

\section{Optimal local approximation spaces and nearly exponential upper
bounds on their accuracy}

\label{sectexponentoverview}

We identify the optimal local approximation spaces for GFEM and provide an
upper bound on the accuracy of their approximation.  We will
consider the case when the exact solution $u_0$ solves the Neumann problem with $f=0$ in \eqref{1.1}. During the course of the exposition we will indicate the modifications needed to treat the corresponding Dirichlet problem. The generalization to the case $f\neq 0$ will be addressed as part of the implementation given in section \ref{implementation}.  

We recall the local approximation spaces introduced in section 2 denoted by $%
\Psi _{i}$ defined on the sets $\omega _{i}$. To fix ideas we will assume
that $\omega _{i}$ are the cubes of a given side length surrounded by a
larger cube $\omega _{i}^{\ast}$. We will distinguish two cases depending on
if the set $\omega _{i},$ lies within the interior of $\Omega$ or if $%
\overline{\omega}_i\cap\partial\Omega\not=\emptyset$. It will be shown that the
overall approach to constructing optimal local approximation spaces for
these two cases is the same. We drop subscripts and consider concentric
cubes $\omega\subset\omega^{\ast}$ with side lengths given by $\sigma$ and $%
\sigma^*=(1+\rho)\sigma$ respectively. In order to introduce the ideas we
suppose first that $\omega$ lies in the interior of $\Omega$ so that $%
\omega\subset\omega^*\subset\Omega$. The energy inner products associated
with these subsets are defined by 
\begin{eqnarray}
(u,v)_{\mathcal{E}(\omega)}=\int_\omega A\nabla u\cdot\nabla v \, dx&&(u,v)_{%
\mathcal{E}(\omega^*)}=\int_{\omega^*} A\nabla u \cdot\nabla v \, dx.
\label{energy}
\end{eqnarray}

We shall utilize $\omega^*$ to construct a finite dimensional approximation
space over $\omega$. For any open subset $S$ of the computational domain $%
\Omega$ we introduce the space of functions $H_A(S)$ defined to be the
functions in $H^1(S)$ that are $A$-harmonic on $S$, i.e., 
$v\in H^{1}(S )$ and 
\begin{equation*}
B(v,\varphi )=0,\,\,\,\,\forall \varphi \in C_{0}^{\infty }(S ).
\end{equation*}
Here $%
H_A(\omega)$ and $H_A(\omega^*)$ contain local information on the
heterogeneities and will be used in the construction of the optimal
local basis. We introduce the quotient of $H_A(\omega^*)$ with respect to the
constant functions denoted by $H_A(\omega^*)/\mathbb{R}$. It is
clear that the solution $u_0$ lies in this local space modulo a constant.

In this method we choose to approximate  elements in the
space of functions $H_A(\omega^*)/\mathbb{R}$ restricted to $\omega$.
This choice is motivated by the Caccioppoli inequality \eqref{3.1}, proved
in Appendex A, which is used to estimate the energy norm over $\omega$ in
terms of the $L^2$ norm over $\omega^*$.  Let $%
P:H_A(\omega^*)/\mathbb{R}\rightarrow H_A(\omega)/\mathbb{R}$ be the restriction
operator such that $P(u)(x)=u(x)$ for all $x\in\omega$ and $u\in
H_A(\omega^*)/\mathbb{R}$. The operator $P$ is compact, this follows from
the Caccioppoli inequality, and the compactness proof is given in the
Appendix, see Theorem \ref{compact}.

Now we approximate by ``$n$'' dimensional subspaces $S(n)\subset
H_A^1(\omega)/\mathbb{R}$. The
accuracy of a particular increasing sequence $\left\{S(n)\right\}_{n=1}^%
\infty$ of local approximation spaces is measured by 
\begin{eqnarray}
d(S(n),\omega)=\sup_{u\in H_A(\omega^*)/\mathbb{R}}\inf_{\chi\in S(n)}\frac{%
\Vert Pu-\chi\Vert_{\mathcal{E}(\omega)}}{\Vert u\Vert_{\mathcal{E}%
(\omega^*)}}.  \label{napprox}
\end{eqnarray}
A sequence of approximation spaces $\hat{S}(n)$ is said to be optimal if it
has an accuracy $d(\hat{S}(n),\omega)$ that satisfies $d(\hat{S}%
(n),\omega)\leq d(S(n),\omega)$, $n=1,2,\ldots$, when compared to any other
sequence of approximation spaces $S(n)$. The problem of finding the family
of optimal local approximation spaces is formulated as follows. Let 
\begin{eqnarray}
d_n(\omega,\omega^*)=\inf_{S(n)}\sup_{u\in H_A(\omega^*)/\mathbb{R}%
}\inf_{\chi\in S(n)}\frac{\Vert Pu-\chi\Vert_{\mathcal{E}(\omega)}}{\Vert
u\Vert_{\mathcal{E}(\omega^*)}}.  \label{nwidth}
\end{eqnarray}
Then the optimal family of approximation spaces $\{\Psi_n(\omega)\}_{n=1}^%
\infty$ satisfy 
\begin{eqnarray}
d_n(\omega,\omega^*)=\sup_{u\in H_A(\omega^*)/\mathbb{R}}\inf_{\chi\in
\Psi_n(\omega)}\frac{\Vert Pu-\chi\Vert_{\mathcal{E}(\omega)}}{\Vert u\Vert_{%
\mathcal{E}(\omega^*)}}.  \label{nwidthopt}
\end{eqnarray}
The quantity $d_n(\omega,\omega^*)$ is known as the Kolomogorov n-width of
the compact operator $P$ see, \cite{Pincus}.

The optimal local approximation space $\Psi_n(\omega)$ for GFEM follows from
general considerations. We introduce the adjoint operator $%
P^*:H_A(\omega)/\mathbb{R}\rightarrow H_A(\omega^*)/\mathbb{R}$ and the operator $P^*P$
is a compact, self adjoint, non-negative operator mapping $H_A(\omega^*)/%
\mathbb{R}$ into itself. We denote the eigenfunctions and eigenvalues of the
problem 
\begin{eqnarray}
P^* P u=\lambda u  \label{eigen}
\end{eqnarray}
by $\{\varphi_i\}$ and $\{\lambda_i\}$ and apply Theorem 2.2, Chapter 4 of 
\cite{Pincus} to see that the optimal subspace $\Psi_n$ is given by the
following theorem.

\begin{theorem}
The optimal approximation space is given by $\Psi_n(\omega)=span\{\psi_1,%
\ldots,\psi_n\}$, where $\psi_i=P\varphi_i$ and $d_n(\omega,\omega^*)=\sqrt{%
\lambda_{n+1}}$. \label{thlambdanplusone}
\end{theorem}

For the case considered here the definitions of $P$ and $P^*$ show that the
optimal subspace and eigenvalues are given by the following explicit
eigenvalue problem.

\begin{theorem}
The optimal approximation space is given by $\Psi_n(\omega)=span\{\psi_1,%
\ldots,\psi_n\}$ where $\psi_i=P\varphi_i$ and $\varphi_i$ and $\lambda_i$
are the first $n$ eigenfunctions and eigenvalues that satisfy 
\begin{eqnarray}
(\varphi_i,\delta)_{\scriptscriptstyle{\mathcal{E}(\omega)}}&=&\lambda_i
(\varphi_i,\delta)_{\scriptscriptstyle{\mathcal{E}(\omega^*)}}%
\hbox{,     
$\forall \delta\in H_A(\omega^*)/\mathbb{R}$}.  \label{concreteevalue}
\end{eqnarray}
\label{concreteinteriorlocal}
\end{theorem}

\noindent \emph{Proof.} The eigenvalue problem \eqref{eigen} is given by 
\begin{eqnarray}
(P^*Pu,\delta)_{\scriptscriptstyle{\mathcal{E}(\omega^*)}}&=&\lambda
(u,\delta)_{\scriptscriptstyle{\mathcal{E}(\omega^*)}}%
\hbox{,      $\forall
\delta\in H_A(\omega^*)/\mathbb{R}$}.  \label{vareigen}
\end{eqnarray}
The theorem follows from \eqref{vareigen} noting that 
\begin{eqnarray}
(P^*Pu,\delta)_{\scriptscriptstyle{\mathcal{E}(\omega^*)}}=(Pu,P\delta)_{%
\scriptscriptstyle{\mathcal{E}(\omega)}}=(u,\delta)_{\scriptscriptstyle{%
\mathcal{E}(\omega)}}.  \label{identities}
\end{eqnarray}

The next theorem provides an upper bound on the rate of convergence for the
optimal local approximation.

\begin{theorem}{\rm Exponential convergence for interior approximations.}

\medskip

For $\epsilon>0$ there is an $N_\epsilon>0$ such that for all $n>N_\epsilon$ 
\begin{eqnarray}
d_n(\omega,\omega^*)\leq e^{-n^{\left(\frac{1}{d+1}-\epsilon\right)}}.
\label{uboundexp}
\end{eqnarray}
\label{boundnearlyexpo}
\end{theorem}

\noindent The index $N_\epsilon$ is constructed explicitly in the proof of
Theorem \ref{boundnearlyexpo} given in the next section. Theorem \ref%
{boundnearlyexpo} shows that the asymptotic convergence rate associated with
the optimal approximation space is nearly exponential for the general class
of $L^\infty(\Omega)$ coefficients belonging to $\mathfrak{C}$. In section %
\ref{locbdry} we identify the optimal local approximation space for the case
when $\omega$ touches the boundary of the computational domain. In that
section we show that the convergence rate is also nearly exponential with
respect to the degrees of freedom of the optimal local approximation space.
We collect our results in section \ref{Macroexponential} and recall Theorem %
\ref{globallocal} to obtain the nearly exponential convergence rate for GFEM
stated in Theorem \ref{Nearlyexpofem}.

\subsection{Local approximation on the interior}
\label{locinterior}

In this section we establish Theorem \ref{boundnearlyexpo}. To do this we
construct a family of approximation spaces that exhibit nearly exponential
convergence in the accuracy of approximation. The convergence rate for this
family delivers an upper bound on the convergence rate for the optimal local
approximation space described by Theorem \ref{concreteinteriorlocal}.

For future reference we introduce the decomposition of $H^1(\omega^*)$
given by 
\begin{eqnarray}
H^1(\omega^*)=H_A^0(\omega^*)+ H_0^1(\omega^*)+
\mathbb{R}.
\nonumber
\end{eqnarray}
Here $H_A^0(\omega^*)$ is the subspace of $H_A(\omega^*)$ given by the 
functions with zero average over $\omega^*$. The spaces $H_A^0(\omega^*)$ and $H_0^1(\omega^*)$ are 
orthogonal with respect to the energy inner
product $(\cdot,\cdot)_{\mathcal{E}(\omega^*)}$. The orthogonal projection
from $H^1(\omega^*)$ onto $H_A^0(\omega^*)$ is denoted by $\mathcal{P}^A$.

The construction of the local approximation space is done iteratively. We
start by introducing the the first $n$ non-constant eigenfunctions $v_i \in H^1(\omega^*)$ of the
Neumann eigenvalue problem 
\begin{eqnarray}
( v_i,w)_{\mathcal{E}(\omega^*)}=\lambda_i \int_{\omega^*} v_i w\,dx,\hbox{  $\forall$ $w\in H^1(\omega^*)$}
\nonumber
\end{eqnarray}
posed over $\omega^{\ast}$, $i=1,\ldots, n$. The subspace spanned by these functions is
denoted by $S_n(\omega^*)$. Next we introduce the span of $A$ harmonic
functions given by 
\begin{eqnarray}
W_n(\omega^{\ast})=\hbox{span}\{w_i\in H_A^0(\omega^{\ast})
:\,w_i=v_i,\hbox{ on  }\partial\omega^{\ast},\,i=1,\ldots n\}
\label{omegaast}
\end{eqnarray}
One readily checks that $W_n(\omega^*)=\mathcal{P}^A S_n(\omega^*)$.

We define the family of approximation spaces $\mathcal{F}_{n}(\omega ,\omega
^{\ast })$ given by the restriction of the elements of $W_n(\omega^{\ast})$
to $\omega$. In what follows we first show that $\mathcal{F}_{n}(\omega
,\omega ^{\ast })$ is a family of local approximation spaces with a rate of
convergence on the order of $n^{-1/d}$, for $d=2,3$. To show this we
introduce a suitable version of the Cacciappoli inequality that bounds
functions in the energy norm over any measurable subset $\mathcal{O}%
\subset\omega^\ast$ for which $dist(\partial \mathcal{O},\partial\omega^%
\ast)>\delta>0$ in terms of the $L^2$ norm over $\omega^{\ast}$.

\begin{lemma}
\label{Theorem 3.1} Let $u$ be $A$-harmonic in $\omega^\ast$ and belong to $%
L^2(\omega^\ast)\cap H_{loc}^1(\omega^\ast)$. Then 
\begin{eqnarray}  \label{3.1}
\parallel u\parallel _{\mathcal{E}(\mathcal{O} )}\leq (2(\beta
)^{1/2}/\delta )\parallel u\parallel _{L^{2}(\omega ^{\ast })}.
\end{eqnarray}
where $\beta $ is defined in \eqref{1.3b}.
\end{lemma}
\noindent The proof of Lemma \ref{Theorem 3.1} is given in the Appendix.
Next we introduce the approximation theorem associated with the 
space $W_n(\omega^*)$ given by

\begin{lemma}
\label{Theorem3.4b} Let $u\in H_A^0(\omega^\ast)$ then there exists
a $v_u\in W_n(\omega^*)$ such that 
\begin{eqnarray}  \label{3.3b}
\Vert u-v_u\Vert_{L^2(\omega^*)}=\inf_{v\in W_{n}(\omega ^{\ast})}\parallel
u-v\parallel _{L^{2}(\omega ^{\ast})}\leq C_n\sigma^{\ast }\frac{
\gamma_d^{1/d}}{\sqrt{4\pi}} \alpha ^{-1/2}\parallel u\parallel _{\mathcal{E}%
(\omega ^{\ast })}
\end{eqnarray}
where $\sigma ^{\ast }$ $\ $is the side length of the cube $\omega ^{\ast },$
$\gamma_d$ is the volume of the unit ball in $\mathbb{R}^d$ and 
$C_{n}= n^{-1/d}(1+o(1))$, for $d=2,3$.
\end{lemma}

\noindent\emph{Proof.} The lemma follows immediately from an upper bound
on the quotient 
\begin{eqnarray}
R=\sup_{u\in H_A^0(\omega^\ast)}\inf_{w\in W_n(\omega^\ast)}\frac{%
\Vert u-w\Vert_{L^2(\omega^\ast)}}{\Vert u\Vert_{\mathcal{E}(\omega^\ast)}}.
\label{quotient}
\end{eqnarray}
Fix $u\in H_A(\omega^\ast)^0$ and denote the projection of $u$ onto 
$W_n(\omega^\ast)$ with respect to the energy norm $\Vert\cdot\Vert_{%
\mathcal{E}(\omega^\ast)}$ by $\mathcal{P}^\mathcal{E} u$. Choosing $w=%
\mathcal{P}^\mathcal{E} u$ and noting that $\Vert(I-\mathcal{P}^\mathcal{E})
u\Vert_{\mathcal{E}(\omega^\ast)}\leq\Vert u\Vert_{\mathcal{E}(\omega^\ast)}$
gives the upper bound 
\begin{eqnarray}  \label{rupper}
R\leq\sup_{u\in H_A^0(\omega^\ast)\perp\,W_n(\omega^\ast)}\frac{%
\Vert u\Vert_{L^2(\omega^\ast)}}{\Vert u\Vert_{\mathcal{E}(\omega^\ast)}}.
\end{eqnarray}
Since $W_n(\omega^\ast)=\mathcal{P}^A S_n(\omega^\ast)$ it follows that 
\begin{eqnarray}
\{u\in H_A^0(\omega^*)\perp\mathcal{P}^A S_n(\omega^*)\}=\{u\in
H_A^0(\omega^*)\perp S_n(\omega^*)\},  \nonumber \\
\{u\in H_A^0(\omega^*)\perp S_n(\omega^*)\}\subset\{u\in
H^1(\omega^*)\perp (S_n(\omega^*)+\mathbb{R})\},
\label{subspaces}
\end{eqnarray}
where the $\perp$ in the second line of \eqref{subspaces} is with respect to the $L^2(\omega^*)$ inner product.
Hence
\begin{eqnarray}  \label{rupper2}
R\leq\sup_{u\in H_A^0(\omega^\ast)\perp\,S_n(\omega^\ast)}\frac{%
\Vert u\Vert_{L^2(\omega^\ast)}}{\Vert u\Vert_{\mathcal{E}(\omega^\ast)}}%
\leq\sup_{u\in H^1(\omega^\ast)\perp\,(S_n(\omega^\ast)+\mathbb{R})}\frac{%
\Vert u\Vert_{L^2(\omega^\ast)}}{\Vert u\Vert_{\mathcal{E}(\omega^\ast)}}=%
\frac{1}{\sqrt{\mu_{n+1}}},
\end{eqnarray}
where $\mu_{n+1}$ is the largest Neumann eigenvalue associated with $%
S_{n+1}(\omega^\ast)$. One has the elementary lower bound $%
\mu_{n+1}\geq\alpha \nu_{n+1}$ where $\nu_{n+1}=4\pi
C_n^{-2}(\sigma^*\gamma_d^{1/d})^{-2}$ is the corresponding Neumann eigenvalue
for the Laplacian with on squares ($d=2$) or cubes ($%
d=3 $). The required upper bound on $R$ now follows and the theorem is
proved.

Now we apply Theorem \ref{Theorem 3.1} to $u-v_u$ on $\omega\subset\omega^{%
\ast}$ and combine it with Theorem \ref{Theorem3.4b} to obtain the following
convergence rate associated with the family of approximation spaces $%
\mathcal{F}_{n}(\omega ,\omega ^{\ast })$ given by

\begin{theorem}
\label{Theorem 3.77} Let $u\in {H}_A^0(\omega ^{\ast })$ then there
exists an approximation $v_u\in \mathcal{F}_n(\omega,\omega ^{\ast })$ for
which 
\begin{eqnarray}  \label{3.74}
\parallel u-v_{u}\parallel _{\mathcal{E}(\omega )}\,=\inf_{v\in \mathcal{F}%
_n(\omega,\omega ^{\ast })}\parallel u-v\parallel _{\mathcal{E}(\omega )}\,
\leq I(\omega ,\omega ^{\ast })\,C_n \parallel u\parallel _{\mathcal{E}%
(\omega ^{\ast })}
\end{eqnarray}
where 
\begin{eqnarray}  \label{3.75}
I(\omega ,\omega ^{\ast })=2\,\frac{\gamma_d^{1/d}}{\sqrt{\pi}}\,\frac{1+\rho%
}{\rho }(\beta /\alpha )^{1/2}&\hbox{ and }& C_n=n^{-1/d}(1+o(1)),\,\,d=2,3,
\end{eqnarray}
where $\gamma_d$ is the volume of the unit ball in dimension $d$.
\end{theorem}

Next we proceed iteratively to construct a family of local approximation
spaces with a rate of convergence that is nearly exponential. For any pair
of two concentric cubes $Q\subset Q^\ast$ we define $\mathcal{F}_{n}(Q
,Q^{\ast })$ to be the space given by the restriction of $W_n(Q^{\ast })$ on 
$Q$. We suppose that $\omega^\ast$ is of side length $\sigma^*$. Let $N>1$
be an integer and we suppose that $\omega$ is of side length $
\sigma$ and $\sigma^*=\sigma(1+\rho)$. Choose $\omega _{j}$, $%
j=1,2,...N$ to be the nested family of concentric cubes with side length $%
\frac{\sigma}{2}(1+\rho(N+1-j)/N)$ for which $\omega=\omega_{N+1}\subset\omega_N\subset%
\omega_{N-1}\subset\cdots\subset\omega_1 = \omega^\ast$. We introduce
the local spaces, $\mathcal{F} _{n}(\omega ,\omega _{N})$, $\mathcal{F}%
_{n}(\omega ,\omega _{N-1})$,...,$\mathcal{F}_{n}(\omega ,\omega _1)$. Put $%
m=N\times n$ and we define the approximation space given by 
\begin{eqnarray}
\mathcal{T}(m,\omega ,\omega ^{\ast })=\mathcal{F}_{n}(\omega ,\omega _{1})+\cdots+ \mathcal{F}%
_{n}(\omega ,\omega _{N}).  \label{iterspacebottumup}
\end{eqnarray}
The convergence rate associated with the local approximation space $\mathcal{%
T}(m,\omega ,\omega ^{\ast })$ is given in the following theorem.

\begin{theorem}
\label{Theorem 3.76} Let $u\in {H}_A^0(\omega ^{\ast })$ and $N\geq
1$ be an integer. Then there exists $z_{u}\in \mathcal{T}(m,\omega ,\omega
^{\ast })$ such that 
\begin{eqnarray}  \label{3.76}
\parallel u-z_{u}\parallel _{\mathcal{E(\omega )}}\leq \varsigma
^{N}\parallel u\parallel _{\mathcal{E}(\omega ^{\ast })}
\end{eqnarray}
and $\varsigma =2\,\frac{\gamma_d^{1/d}}{\sqrt{\pi}}\,\frac{1+\rho}{\rho}
\,N(\beta /\alpha )^{1/2}C_n$.
\end{theorem}

\noindent\emph{Proof.} In what follows we make the identification $%
\omega_0=\omega^*$ and $\omega_{N+1}=\omega$. From Theorem \ref{Theorem 3.77}
we have that there exists $v_1\in\mathcal{F}_{n}(\omega_1 ,\omega^\ast)$
such that 
\begin{eqnarray}  \label{3.7a}
\Vert u-v_1\Vert_{\mathcal{E}(\omega_1)}\leq \varsigma \parallel u\parallel
_{\mathcal{E}(\omega ^{\ast })}
\end{eqnarray}
Suppose next that for $m=1,\ldots, j$ there are functions $v_m\in\mathcal{F}%
_n(\omega_m,\omega_{m-1})$ such that 
\begin{eqnarray}  \label{inductionhyp}
\parallel u-\sum_{m=1}^j v_m\parallel _{\mathcal{E}(\omega_j)}\leq
\varsigma^j \parallel u\parallel _{\mathcal{E}(\omega ^{\ast })}
\end{eqnarray}
Applying Theorem \ref{Theorem 3.77} we see that there exists a $v_{j+1}\in%
\mathcal{F}_{n}(\omega_{j+1} ,\omega_j)$ for which 
\begin{eqnarray}  \label{3.7b}
\Vert u-(\sum_{m=1}^j v_m)-v_{j+1}\Vert_{\mathcal{E}(\omega_{j+1})}\leq
\varsigma \parallel u-\sum_{m=1}^j v_m\parallel _{\mathcal{E}(\omega_j)}
\end{eqnarray}
and the induction step goes through. Choosing $z_{u}=\sum_{m=1}^N v_m$
delivers 
\begin{eqnarray}  \label{3.8}
\parallel u-z_{u}\parallel _{\mathcal{E(\omega )}}\leq \varsigma
^{N}\parallel u\parallel _{\mathcal{E(\omega }^{\ast })}
\end{eqnarray}
and the theorem follows noting that $z_u$ belongs to $\mathcal{T}(m,\omega
,\omega ^{\ast })$.

Next we make a choice for $N$. We choose $N$ to be the largest integer less
than or equal to $n^\gamma$ for $0<\gamma$. Thus $m\leq n^{\gamma+1}$ and $m^{\frac{1}{%
\gamma+1}}\leq n$ and it follows that $n^{-\frac{1}{d}}\leq m^{-\frac{1}{%
d(\gamma+1)}}$, and $N\leq m^\frac{\gamma}{\gamma+1}$. On applying these
inequalities  we obtain
\begin{eqnarray}
\varsigma^N\leq\exp{\left\{-m^{\frac{\gamma}{\gamma+1}}\left(-\ln{K}+\frac{%
1/d-\gamma}{\gamma+1}\ln{m}\right)\right\}}  \label{exponent}
\end{eqnarray}
where $K=2\,\frac{\gamma_d^{1/d}}{\sqrt{\pi}}(\frac{\beta}{\alpha})^{1/2}(
\frac{1+\rho}{\rho})$. It is evident that decay occurs for the choice $0\leq\gamma<%
\frac{1}{d}$ and 
\begin{eqnarray}
\varsigma^N < e^{-m^{\frac{\gamma}{\gamma+1}}}  \label{valueofN}
\end{eqnarray}
for $m>N=(Ke)^{(\gamma+1)/(1/d-\gamma)}$. We set $\ell={\rm dimension}\{\mathcal{T}(m,\omega ,\omega ^{\ast })\}$  
and Theorem \ref{Theorem 3.76} together with
\eqref{valueofN} imply 
\begin{eqnarray}
d_\ell(\omega,\omega^*)\leq\sup_{u\in H_A(\omega^*)/\mathbb{R}}\inf_{\chi\in 
\mathcal{T}(m,\omega ,\omega ^{\ast })}\frac{\Vert u-\chi\Vert_{\mathcal{E}%
(\omega)}}{\Vert u\Vert_{\mathcal{E}(\omega^*)}}\leq e^{-\ell^{\frac{\gamma}{%
\gamma+1}}}.  \label{upbd}
\end{eqnarray}
for $\ell>N$ and Theorem \ref{boundnearlyexpo} is proved.

\subsection{Local approximation at the boundary}

\label{locbdry}

%%%%%%%%%%%%%%%%%%%%%%%%%%%

\begin{figure}[t]
%[ht]
\centering
\scalebox{0.3}{\includegraphics{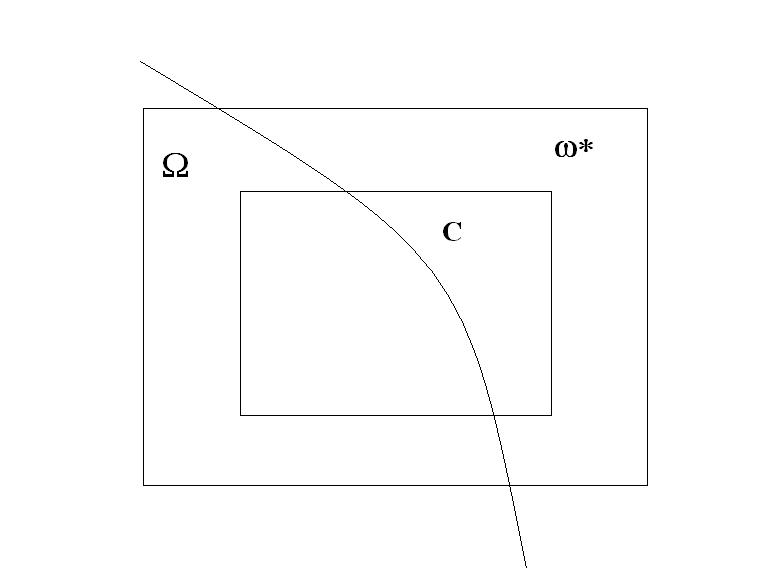}}
\caption{Truncated cube $\protect\omega=C\cap\Omega$.}
\label{Trunc}
\end{figure}

Consider two concentric cubes $C\subset\omega^*$ of side lengths $\sigma$
and $\sigma^*=(1+\rho)\sigma$ respectively. We suppose that 
$\omega^*\cap\Omega\not=\emptyset$ and $\omega^*\cap\partial\Omega\not=\emptyset$.
The truncated cube $\omega$ is defined to
be $\omega=C\cap\Omega$, and $\partial\omega\cap\partial\Omega\not=\emptyset$, see Figure \ref{Trunc}. For this case we will
assume that $\partial\Omega$ is $C^1$, i.e., the boundary can be represented
locally as the graph of a $C^1$ function. The method presented here applies
to both Dirichlet and Neumann boundary value problems. We will illustrate
the ideas for the Neumann problem and make references to the Dirichlet
problem when appropriate. Given a function $u\in H_A(\Omega)$ the goal is to
provide a local approximation to $u$ in $\omega$. To this end we form a
local particular solution $u_p$ given by the $A$-harmonic function that
satisfies $n\cdot A\nabla u_p=g$ on $\omega^*\cap\partial\Omega$ and $u_p=0$ on $\partial\omega^*\cap\Omega$.
Writing $u=u_p+u_0$ we see that $n\cdot A\nabla u_0=0$ on $\omega^*\cap\partial\Omega$ 
and $u_0=u$ on $\partial\omega^*\cap\Omega$. If
instead we have Dirichlet data then the particular solution $u_p$ satisfies 
$u_p=q=u$ on $\omega^*\cap\partial\Omega$ and $u_p=0 $ on $\partial\omega^*\cap\Omega$. The objective of this section is to
find the optimal family of local approximation spaces that give the best
approximation to $u_0=u-u_p$ in the energy norm over the set $\omega$. To
this end we introduce the the space of functions $H_{A,0}(\Omega\cap\omega^%
\ast)$ given by all functions $v$ in $H^1(\Omega\cap\omega^\ast)$ that are $%
A $-harmonic on $\Omega\cap\omega^\ast$ and for which 
$$\partial_\nu v\equiv n\cdot A\nabla v=0$$ on 
$\partial\Omega\cap\omega^{\ast}$. The analogous space of functions defined
on $\omega$ is denoted by $H_{A,0}(\omega)$. Since we approximate functions with respect to the energy norm we introduce the quotient space of $H_{A,0}(\Omega\cap\omega^%
\ast)$ with respect to the constant functions denoted by $H_{A,0}(\Omega\cap\omega^%
\ast)/\mathbb{R}$.

Now we introduce  $P:H_{A,0}(\omega^*\cap%
\Omega)/\mathbb{R}\rightarrow H_{A,0}(\omega)/\mathbb{R}$ given by the restriction operator defined by $%
P(u)(x)=u(x)$ for all $x\in\omega$ and $u\in H_{A,0}(\omega^*\cap\Omega)/\mathbb{R}$.
The operator $P$ is compact, this follows from Lemma \ref{CompactBDRY}
given in the Appendix. Let $S(n)$ be any finite dimensional subspace of $%
H_{A,0}(\omega)/\mathbb{R}$ and the problem of finding the family of optimal local
approximation spaces is formulated in terms of the n-width of $P$. Let 
\begin{eqnarray}
d_n(\omega,\omega^*\cap\Omega)=\inf_{S(n)}\sup_{u\in
H_{A,0}(\omega^*\cap\Omega)/\mathbb{R}}\inf_{\chi\in S(n)}\frac{\Vert Pu-\chi\Vert_{%
\mathcal{E}(\omega)}}{\Vert u\Vert_{\mathcal{E}(\omega^*\cap\Omega)}}.
\label{nwidthBDRY}
\end{eqnarray}
Then the optimal family of boundary approximation spaces $%
\{\Psi_n(\omega)\}_{n=1}^\infty$ for GFEM satisfy 
\begin{eqnarray}
d_n(\omega,\omega^*\cap\Omega)=\sup_{u\in
H_{A,0}(\omega^*\cap\Omega)/\mathbb{R}}\inf_{\chi\in \Psi_n(\omega)}\frac{\Vert
Pu-\chi\Vert_{\mathcal{E}(\omega)}}{\Vert u\Vert_{\mathcal{E}%
(\omega^*\cap\Omega)}}.  \label{nwidthoptBDRY}
\end{eqnarray}
Proceeding as before we introduce the adjoint operator $P^*:H_A(\omega)/\mathbb{R}
\rightarrow H_A(\omega^*\cap\Omega)/\mathbb{R}$ and the operator $P^*P$ is a compact
operator mapping $H_{A,0}(\omega^*\cap\Omega)/\mathbb{R}$ into itself. Similar
arguments show that the optimal approximating spaces are given by the
following theorem.

\begin{theorem}
The optimal approximation space is given by $\Psi_n(\omega)=span\{\psi_1,%
\ldots,\psi_n\}$ where $\psi_i=P\varphi_i$ and $\varphi_i\in
H_{A,0}(\omega^*\cap\Omega)/\mathbb{R}$ and $\lambda_i$ are the first $n$
eigenfunctions and eigenvalues that satisfy 
\begin{eqnarray}
(\varphi_i,\delta)_{\scriptscriptstyle{\mathcal{E}(\omega)}}&=&\lambda_i
(\varphi_i,\delta)_{\scriptscriptstyle{\mathcal{E}(\omega^*\cap\Omega)}}%
\hbox{,      $\forall \delta\in H_{A,0}(\omega^*\cap\Omega)/\mathbb{R}$}.
\label{concreteevalueBDRY}
\end{eqnarray}
\label{concreteinteriorlocalBDRY}
\end{theorem}

The next theorem provides an upper bound on the rate of convergence for the
optimal local boundary approximation.

\begin{theorem}{\rm Exponential convergence at the boundary.}

\medskip

For $\epsilon>0$ there is an $N_\epsilon>0$ such that for all $n>N_\epsilon$ 
\begin{eqnarray}
d_n(\omega,\omega^*\cap\Omega)\leq e^{-n^{\left(\frac{1}{d+1}%
-\epsilon\right)}}.  \label{uboundexpBDRY}
\end{eqnarray}
\label{boundnearlyexpoBDRY}
\end{theorem}

\noindent Theorem \ref{boundnearlyexpoBDRY} shows that the asymptotic
convergence rate associated with the optimal boundary approximation space is
also nearly exponential for the general class of $L^\infty(\omega^*)$
coefficients $\mathfrak{C}$.

We now give the prof of Theorem \ref{boundnearlyexpoBDRY}. 
The subspace of A-harmonic functions defined over $\Omega\cap\omega^*$ with zero mean is denoted by $H_A^0(\Omega\cap\omega^*)$
and the subspace of elements belonging to $H_{A,0}(\Omega\cap\omega^\ast)$ 
with zero mean over $\Omega\cap\omega^*$ is denoted by $H_{A,0}^0(\Omega\cap\omega^\ast)$.
We introduce the $ L^2(\Omega\cap\omega^\ast)$ norm closure 
of $H_{A,0}^0(\Omega\cap\omega^\ast)$
denoted by $\overline{H}_{A,0}^0(\Omega\cap\omega^\ast)$. Useful properties of functions in $
\overline{H}_{A,0}^0(\Omega\cap\omega^\ast)$ are listed below in Theorem \ref{Closed} at 
the end of this section. It is shown there that $v\in\overline{H}%
_{A,0}^0(\Omega\cap\omega^\ast)$ implies that $\partial_\nu v=0$ on $\partial\Omega\cap%
\omega^\ast$ and that $v$ is $A$-harmonic in $\Omega\cap\omega^\ast$.

Next we introduce the the first $n$ non-constant eigenfunctions of the Neumann
eigenvalue problem $\mathrm{div}(A\,\nabla v_i)=-\lambda_i v_i$ over $%
\Omega\cap\omega^{\ast}$, $i=1,\ldots, n$. The subspace spanned by these
functions is denoted by $S_n(\Omega\cap\omega^*)$. Next we introduce the
span of $A$ harmonic functions given by 
\begin{eqnarray}
W_n(\Omega\cap\omega^{\ast})=\hbox{span}\{w_i\in
H_A^0(\Omega\cap\omega^{\ast}):\,w_i=v_i,\hbox{ on 
}\partial\omega^{\ast}\cap\Omega\cup\partial\Omega\cap\omega^*,\,i=1,\ldots n\}.  \label{omegaastOmega}
\end{eqnarray}
For future reference we note the  decomposition of 
$H^1(\Omega\cap\omega^*)$ given by 
$$H^1(\Omega\cap\omega^*)=H_A^0(\Omega\cap%
\omega^*)+H_0^1(\Omega\cap\omega^*)+\mathbb{R}.$$ 
Here the first two subspaces are orthogonal with
respect to the energy inner product $(\cdot,\cdot)_{\mathcal{E}%
(\Omega\cap\omega^*)}$. The orthogonal projection from $H^1(\Omega\cap%
\omega^*)$ onto $H_A^0(\Omega\cap\omega^*)$ is denoted by $\mathcal{P}^A$. It
is easily verified that $W_n(\Omega\cap\omega^*)=\mathcal{P}^A
S_n(\Omega\cap\omega^*)$. 

We now $L^2$ project this space onto $%
\overline{H}_{A,0}^0(\Omega\cap\omega^\ast)$. The projection
operator mapping $L^2(\Omega\cap\omega^*)$ onto $\overline{H}%
_{A,0}^0(\Omega\cap\omega^\ast)$ is denoted by $\mathcal{P}_0$ and 
\begin{eqnarray}
\Vert v-\mathcal{P}_0 v\Vert_{L^2(\omega^\ast\cap\Omega)}= \inf_{w\in%
\overline{H}_{A,0}^0(\Omega\cap\omega^\ast)}\Vert
v-w\Vert_{L^2(\omega^\ast\cap\Omega)}\hbox{  and}  \notag \\
\Vert v\Vert^2_{L^2(\omega^\ast\cap\Omega)}=\Vert \mathcal{P}_0
v\Vert^2_{L^2(\omega^\ast\cap\Omega)}+\Vert (I-\mathcal{P}_0)
v\Vert^2_{L^2(\omega^\ast\cap\Omega)}  \label{argminbdry}
\end{eqnarray}
In what follows the local approximations will be chosen from the local
function space $\mathcal{P}_0 W_n(\omega^{\ast}\cap\Omega)$ restricted to
the set $\omega$.

As before we suppose that the side length of $\omega^\ast$ is given by $%
\sigma^*=\sigma(1+\rho)$ and $\omega=C\cap\Omega$ where $C$ is the
concentric sub-cube of $\omega^*$ and the side length of $C$ is $\sigma$. We
let $\omega^{\ast\ast}=\omega^\ast\cap\Omega$ and $\sigma\rho/2=dist(%
\partial\omega\cap\Omega,\partial\omega^{\ast\ast}\cap\Omega)>0$. From the
smoothness assumption on $\partial \Omega $ there exists a dimorphism ${%
\mathbf{x}=\mathbf{\Psi}}(\mathbf{y})$ of class $C^{1}$ for which $\mathbf{%
\Psi}^{-1}$maps $\omega ^{\ast\ast }$ onto $\widetilde{\omega ^{\ast\ast }}$
with $\partial \widetilde{\omega ^{\ast\ast }}\cap \partial \widetilde{%
\Omega }$ being a part of the plane $y_{n}=0$. We extend $\widetilde{\omega
^{\ast\ast }}$ across $y_n=0$ by reflection. We denote the image of this
extension under ${\mathbf{\Psi}}$ by $\omega_E^{\ast\ast}$ and define $%
\omega_E^\ast$ to be given by the union $\omega^{\ast\ast}\cup(\partial%
\omega^{\ast\ast}\cap\partial\Omega)\cup\omega_E^{\ast\ast}$.

In what follows we extend $u-\mathcal{P}_0 v$ across the boundary $%
\partial\Omega$ as an $A$-Harmonic function over $\omega_E^\ast$ and apply
the Caciappoli inequality to recover the following theorem.

\begin{lemma}
{\rm $L^2$ projection of local fields at the boundary.}\newline
Suppose $\omega$ is a truncated cube with part of its boundary given by $%
\partial\Omega$ and suppose there exists a dimorphism $\mathbf{x}=\mathbf{%
\Psi}(\mathbf{y})$ of class $C^{1}$ for which $\mathbf{\Psi}^{-1}$ maps $%
\omega ^{\ast\ast }$ onto $\widetilde{\omega ^{\ast\ast }}$ with $\partial 
\widetilde{\omega ^{\ast\ast }}\cap \partial \widetilde{\Omega }$ being a
part of the plane $y_{n}=0$. Let $u\in H_{A,0}(\Omega\cap \omega^\ast)$,
then there is a constant $C>0$ depending only on $\partial\Omega$ such that
given $v\in L^2(\omega^\ast\cap\Omega)$ the projection $\mathcal{P}_0 v$
onto $\overline{H}_{A,0}^0(\Omega\cap\omega^\ast)$ satisfies 
\begin{eqnarray}
\parallel u-\mathcal{P}_0 v \parallel _{\mathcal{E}(\omega)}\leq
4C(\beta)^{1/2}\frac{1}{\sigma\rho} \parallel u- v\parallel _{L^{2}(\omega
^{\ast }\cap \Omega )}.  \label{4.23}
\end{eqnarray}
\label{Theorem4.4}
\end{lemma}

\emph{Proof.} We set $v^\ast=\mathcal{P}_0 v$ and and observe from Theorem %
\ref{Closed} that $\partial_\nu (u-v^\ast)$ vanishes on $\partial\Omega\cap\omega^\ast$.
Let $\widetilde{u}=u(\mathbf{\Psi}(\mathbf{y}))$ and $\widetilde{v^\ast}%
=v^\ast(\mathbf{\Psi}(\mathbf{y}))$ and note that $(\widetilde{u}-\widetilde{%
v^\ast})$ is an $\widetilde{A}-$ harmonic function
in $\widetilde{\omega^{\ast\ast}}$ where $\widetilde{A}(\mathbf{y})=[\nabla%
\mathbf{\Psi}]^{-1}A(\mathbf{\Psi}(\mathbf{y}))[\nabla\mathbf{\Psi}]^{-1}$ and $\tilde{n}\cdot\tilde{A}(\tilde{u}-\widetilde{v^\ast})=0$. Here $\tilde{n}$ is the unit normal to $y_n=0$
pointing into $y_n<0$.
Since $\tilde{n}\cdot\tilde{A}(\tilde{u}-\widetilde{v^\ast})=0$ we apply standard arguments to extend $\widetilde{A}$ across $y_n=0$  so that $(\widetilde{u}-\widetilde{v^\ast})$ is extended 
across $y_n=0$ outside $\widetilde{\Omega }$ as an $\widetilde{A}-$harmonic
function. We set $y'=(y_1,\ldots,y_{n-1})$ and write $y=(y',y_n)$. For $y_n<0$ we extend $\widetilde{A}$ across $y_n=0$ into $y_n<0$ according to: 1) $\widetilde{A}_{ij}(y',-y_n)$, for all $i=j$, $j=1,\ldots,n$, 2) $\widetilde{A}_{ij}(y',-y_n)$ for all $j\not=n,i\not=n$,   3)  $-\widetilde{A}_{ij}(y',-y_n)$ for $j=n$, $i<n$,  and 4) $-\widetilde{A}_{ij}(y',-y_n)$ for $i=n$ and all $j<n$. 

We map back to obtain an extension of  $A$ across $\partial\Omega$ and recover an $%
A$-harmonic extension of the function $u-v^\ast$ on $\omega_E^\ast$. From
Theorem \ref{Closed} it follows that $u-v^\ast\in L^2(\omega_E^\ast)\cap
H^1_{loc}(\omega_E^\ast)$ and we apply Theorem \ref{Theorem 3.1} to ${\omega}%
\subset\omega_E^\ast$, to find that 
\begin{eqnarray}
&&\parallel u-v^\ast\parallel _{\mathcal{E(\omega)}}\leq 4(\beta)^{1/2}\frac{%
1}{\sigma\rho} \parallel u-v^\ast\parallel _{L^{2}(\omega_E ^{\ast })} 
\notag \\
&&\leq 4C(\beta)^{1/2}\frac{1}{\sigma\rho} \parallel u-v^\ast\parallel
_{L^{2}(\omega ^{\ast }\cap\Omega)}  \notag \\
&&=4C(\beta)^{1/2}\frac{1}{\sigma\rho} \parallel \mathcal{P}_0(u-v)\parallel
_{L^{2}(\omega ^{\ast }\cap\Omega)}\leq 4C(\beta)^{1/2}\frac{1}{\sigma\rho}
\parallel u-v\parallel _{L^{2}(\omega ^{\ast }\cap\Omega)}
\label{outeromega}
\end{eqnarray}
and the theorem is proved. We point out that the analogous theorem holds for the Dirichlet problem
and can be proved using similar arguments. 

In what follows it is always assumed that $\partial\Omega\cap\omega^{\ast}$
can be flattened according to the hypothesis of Lemma \ref{Theorem4.4}.
Next we introduce the approximation theorem associated with the space $%
\mathcal{P}_0 W_n(\Omega\cap\omega^*)$ given by

\begin{lemma}
\label{Theorem3.41b} Let $u\in H_{A,0}^0(\Omega\cap\omega^\ast)$ then there
exists a $v_u\in \mathcal{P}_0 W_n(\Omega\cap\omega^*)$ such that 
\begin{eqnarray}  \label{3.31b}
\Vert u-v_u\Vert_{L^2(\Omega\cap\omega^*)}=\inf_{v\in \mathcal{P}_0
W_{n}(\Omega\cap\omega ^{\ast})}\parallel u-v\parallel
_{L^{2}(\Omega\cap\omega ^{\ast})}\leq C_n\sigma^{\ast }\frac{\gamma_d^{1/d}%
}{\sqrt{4\pi}} \alpha ^{-1/2}\parallel u\parallel _{\mathcal{E}%
(\Omega\cap\omega ^{\ast })}
\end{eqnarray}
where $\sigma ^{\ast }$is the side length of the cube $\omega ^{\ast },$ $%
\gamma_d$ is the volume of the unit ball in $\mathbb{R}^d$ and $C_{n}=
n^{-1/d}(1+o(1))$, for $d=2,3$.
\end{lemma}

\emph{Proof.} The theorem follows immediately from an upper bound on the
quotient 
\begin{eqnarray}
R=\sup_{u\in H_{A,0}^0(\Omega\cap\omega^\ast)}\inf_{w\in \mathcal{P}_0
W_n(\Omega\cap\omega^\ast)}\frac{\Vert u-w\Vert_{L^2(\Omega\cap\omega^\ast)}%
}{\Vert u\Vert_{\mathcal{E}(\Omega\cap\omega^\ast)}}.  \label{quotient1}
\end{eqnarray}
%For $w\in\mathcal{P}_0 W_n(\Omega\cap\omega^*)$ we write it as $w=\mathcal{P}_0 g$, for $g\in W_n(\Omega\cap\omega^*)$. 
Fix $u\in H_{A,0}^0(\Omega\cap\omega^\ast)$ and for every $w\in\mathcal{P}_0
W_n(\Omega\cap\omega^*)$ one has $g\in W_n(\Omega\cap\omega^*)$ such that $w=%
\mathcal{P}_0 g$ and 
\begin{eqnarray}
\inf_{w\in \mathcal{P}_0
W_n(\Omega\cap\omega^\ast)}\Vert u-w\Vert_{L^2(\Omega\cap\omega^\ast)}&=&\inf_{g\in
W_n(\Omega\cap\omega^\ast)}\Vert \mathcal{P}%
_0(u-g)\Vert_{L^2(\Omega\cap\omega^\ast)}  \notag \\
&\leq&\inf_{g\in
W_n(\Omega\cap\omega^\ast)}\Vert u-g\Vert_{L^2(\Omega\cap\omega^\ast)}.  \label{projineq}
\end{eqnarray}
Thus
\begin{eqnarray}
R\leq\sup_{u\in H_{A,0}^0(\Omega\cap\omega^\ast)}\inf_{g\in
W_n(\Omega\cap\omega^\ast)}\frac{\Vert u-g\Vert_{L^2(\Omega\cap\omega^\ast)}%
}{\Vert u\Vert_{\mathcal{E}(\Omega\cap\omega^\ast)}}.  \label{quotient2}
\end{eqnarray}
Denote the projection of $u$ onto $W_n(\omega^\ast)$ with respect to the
energy norm $\Vert\cdot\Vert_{\mathcal{E}(\Omega\cap\omega^\ast)}$ by $%
\mathcal{P}^\mathcal{E} u$. Choosing $g=\mathcal{P}^\mathcal{E} u$ and
noting that $\Vert(I-\mathcal{P}^\mathcal{E}) u\Vert_{\mathcal{E}%
(\Omega\cap\omega^\ast)}\leq\Vert u\Vert_{\mathcal{E}(\Omega\cap\omega^%
\ast)} $ gives the upper bound 
\begin{eqnarray}  \label{rupper7}
R\leq\sup_{u\in
H_{A,0}^0(\Omega\cap\omega^\ast)\perp\,W_n(\Omega\cap\omega^\ast)}\frac{\Vert
u\Vert_{L^2(\Omega\cap\omega^\ast)}}{\Vert u\Vert_{\mathcal{E}%
(\Omega\cap\omega^\ast)}}.
\end{eqnarray}
Now  $H_{A,0}^0(\Omega\cap\omega^\ast)\subset H_A^0(\Omega\cap\omega^\ast)$ so 
\begin{eqnarray}  \label{rupper3}
R\leq\sup_{u\in (H_A^0(\Omega\cap\omega^\ast))\perp\,W_n(\Omega%
\cap\omega^\ast)}\frac{\Vert u\Vert_{L^2(\Omega\cap\omega^\ast)}}{\Vert
u\Vert_{\mathcal{E}(\Omega\cap\omega^\ast)}}.
\end{eqnarray}
Since $W_n(\Omega\cap\omega^\ast)=\mathcal{P}^A S_n(\Omega\cap\omega^\ast)$
it follows that 
\begin{eqnarray}
\{u\in H_A^0 (\Omega\cap\omega^*)\perp\mathcal{P}^A
S_n(\Omega\cap\omega^*)\}&=&\{u\in H_A^0(\Omega\cap\omega^*)\perp
S_n(\Omega\cap\omega^*)\}  \notag \\
&&\subset\{u\in H^1(\Omega\cap\omega^*)\perp (S_n(\Omega\cap\omega^*)+\mathbb{R})\}.
\label{subspaces3}
\end{eqnarray}
Here on the second line of \eqref{subspaces3} the $\perp$ is with respect to the $L^2(\Omega\cap\omega^*)$ inner product. It now follows that 
\begin{eqnarray}  \label{rupper5}
R&\leq&\sup_{u\in H_A^0(\Omega\cap\omega^\ast)\perp\,S_n(\Omega\cap%
\omega^\ast)}\frac{\Vert u\Vert_{L^2(\Omega\cap\omega^\ast)}}{\Vert u\Vert_{%
\mathcal{E}(\Omega\cap\omega^\ast)}}  \notag \\
&\leq&\sup_{u\in H^1(\Omega\cap\omega^\ast)\perp\,(S_n(\Omega\cap%
\omega^\ast)+\mathbb{R})}\frac{\Vert u\Vert_{L^2(\Omega\cap\omega^\ast)}}{\Vert u\Vert_{%
\mathcal{E}(\Omega\cap\omega^\ast)}}=\frac{1}{\sqrt{\mu_{n+1}}},
\end{eqnarray}
where $\mu_{n+1}$ is the largest Neumann eigenvalue associated with $%
S_{n+1}(\omega^\ast)$. One has the elementary lower bound $%
\mu_{n+1}\geq\alpha \nu_{n+1}$, where $\nu_{n+1}$ is the associated Neumann eigenvalue for the Laplacian on $\omega^*\cap\Omega$. 
For this case Weyl's theorem \cite{Weyl}
gives 
\begin{eqnarray}
\nu_{n+1}= 4\pi\left(\frac{n}{\gamma_d|\Omega\cap\omega^*|}
\right)^{2/d}+o(n^{2/d}).  \label{weyl}
\end{eqnarray}
The upper bound on $R$ now follows from \eqref{weyl} together with the
inequality $|\Omega\cap\omega^*|\leq(\sigma^*)^d$ and the theorem is proved.

We introduce the local space near the boundary given by $\partial\mathcal{F}%
_{n}(\omega ^{\ast }\cap\Omega)=\mathcal{P}_0 W_{n}(\omega^\ast\cap\Omega)$
and define the local approximation space $\partial\mathcal{F}_{n}(\omega
,\omega ^{\ast }\cap\Omega)$ to be given by the restriction of $\partial%
\mathcal{F}_{m}(\omega ^{\ast }\cap\Omega)$ on $\omega $.

Now we apply Lemma \ref{Theorem4.4} to $u-v_u$ on $\omega\subset\Omega\cap%
\omega^{\ast}$ and combine it with Lemma \ref{Theorem3.41b} to obtain the
following convergence rate associated with the family of approximation
spaces $\partial\mathcal{F}_{n}(\omega ,\omega ^{\ast }\cap\Omega)$.

\begin{lemma}
\label{Theorem 4.77} Let $u\in H_{A,0}^0(\omega^{\ast }\cap\Omega)$, then
there exists an approximation $v_u\in \partial\mathcal{F}_n(\omega,\omega
^{\ast }\cap\Omega)$ for which 
\begin{eqnarray}  \label{4.174}
\parallel u-v_{u}\parallel _{\mathcal{E}(\omega )}\,=\inf_{w\in \partial%
\mathcal{F}_n(\omega,\omega ^{\ast }\cap\Omega)}\parallel u-w\parallel _{%
\mathcal{E}(\omega )}\, \leq I(\omega ,\omega ^{\ast }\cap\Omega)\,C_n
\parallel u\parallel _{\mathcal{E}(\omega ^{\ast }\cap\Omega)}
\end{eqnarray}
where 
\begin{eqnarray}  \label{4.175}
I(\omega ,\omega ^{\ast }\cap\Omega)=8C\,\frac{\gamma_d^{1/d}}{\sqrt{\pi}}\,%
\frac{1+\rho}{\rho }(\beta /\alpha )^{1/2}&\hbox{ and }& C_n=n^{-1/d}(1+o(1)),\,\,d=2,3,
\end{eqnarray}
where $\gamma_d$ is the volume of the unit ball in dimension $d$ and $C$
depends only upon $\partial\Omega$.
\end{lemma}

Now we proceed iteratively to construct a family of local approximation
spaces with a rate of convergence that is nearly exponential. For any pair
of two concentric cubes $Q\subset \tilde{Q}$ such that their intersections $%
\omega=Q\cap\Omega$ and $\tilde{\omega}=\tilde{Q}\cap\Omega$ have nonzero
volume we define $\partial\mathcal{F}_{n}(\omega ,\tilde{\omega})$ to be the
space given by the restriction of $\mathcal{P}_0 W_n(\tilde{\omega})$ on $%
\omega$. We recall that the two concentric cubes $C\subset \omega^\ast$ are of side 
length $\sigma$ and $\sigma^*=\sigma(1+\rho)$ respectively
and $\omega=C\cap\Omega$, see Figure \ref{Trunc}. Let $%
N>1 $ be an integer and consider the nested family of concentric cubes
$Q _{j}$, $j=1,2,...N+1$ with $Q_{j+1}\subset Q_{j}$ and $Q_1=\omega^*$ and $Q_{N+1}=C$.
The side lengths of $Q_j$ are given by $\sigma(1+\rho(N+1-j)/N))/2$. 
Set $\omega_j=Q_j\cap\Omega$ to obtain $%
\omega=\omega_{N+1}\subset\omega_N\subset\cdots\subset\omega_1=\omega^\ast\cap\Omega$. 
We introduce the local spaces, $\partial\mathcal{F}
_{n}(\omega ,\omega _{N})$, $\partial\mathcal{F}_{n}(\omega ,\omega _{N-1})$%
,...,$\partial\mathcal{F}_{n}(\omega ,\omega _1)$. Put $m=N\times n$ and we
define the approximation space 
\begin{eqnarray}
\Psi=\partial\mathcal{T}(m,\omega ,\omega ^{\ast }\cap\Omega)=\partial\mathcal{F}%
_{n}(\omega ,\omega _{1})+\cdots+ \partial\mathcal{F}_{n}(\omega
,\omega _{N}).  \label{iterspacebottumupbdry}
\end{eqnarray}
The convergence rate associated with the local approximation space $\partial%
\mathcal{T}(m,\omega ,\omega ^{\ast }\cap\Omega)$ is given in the following
theorem.

\begin{theorem}
\label{Theorem 4.876} Let $u\in H_{A,0}^0(\omega ^{\ast }\cap\Omega)$, then
there exists $z_{u}\in \Psi=\partial\mathcal{T}(m,\omega ,\omega ^{\ast
}\cap\Omega)$ such that 
\begin{eqnarray}  \label{4.276}
\parallel u-z_{u}\parallel _{\mathcal{E(\omega )}}\leq \varsigma
^{N}\parallel u\parallel _{\mathcal{E}(\omega ^{\ast }\cap\Omega)}
\end{eqnarray}
and $\varsigma =8\,C\,\frac{\gamma_d^{1/d}}{\sqrt{\pi}}\frac{1+\rho}{\rho}
N(\beta /\alpha )^{1/2}C_n$.
\end{theorem}

\emph{Proof.} The proof is by induction and is identical to the proof of
Theorem \ref{Theorem 3.76}.

Theorem \ref{boundnearlyexpoBDRY} now follows on choosing the appropriate $N$
and using arguments identical to those used to establish Theorem \ref%
{boundnearlyexpo}.

We conclude by stating and proving the following theorem.

\begin{theorem}
The set $\overline{H}_{A,0}^0(\omega ^{\ast }\cap\Omega)$ is a subspace of the
space of $A$-harmonic functions belonging to $L^2(\Omega\cap\omega^\ast)$.
Functions $v$ belonging to $\overline{H}_{A,0}(\omega ^{\ast }\cap\Omega)$
have the following local properties. For any open subset $\mathcal{O}%
\subset\omega^*\cap\Omega$,

\begin{enumerate}
\item $v\in H^1(\mathcal{O})$ for any $\mathcal{O}\subset\omega^*\cap\Omega$%
, such that $dist(\partial\mathcal{O},\partial\omega^*\cap\Omega)>0$, and

\item if $\partial\mathcal{O}\cap(\partial\Omega\cap\omega^*)\not=\emptyset$
then $\partial_\nu v=0$ on $\partial\mathcal{O}\cap(\partial\Omega\cap\omega^*)$.
\end{enumerate}

\label{Closed}
\end{theorem}

\emph{Proof.} 
%Denote the the space of $A$-harmonic functions $v$ belonging to $L^2(\Omega\cap\omega^\ast)\cap H^1_{loc}(\Omega\cap\omega^\ast)$ with $v=0$ on $\partial\Omega\cap\omega^\ast$ by $\mathcal{K}_0$.
%Now consider $u_n\in{H}_{A,0}(\omega }^{\ast }\cap\Omega)$ for which $u_n\rightarrow u_\infty$ in $L^2(\omega }^{\ast }\cap\Omega)$.
Given $u_\infty\in\overline{H}_{A,0}^0(\omega^{\ast }\cap\Omega)$ then there is a sequence $u_n\in H_{A,0}^0(\omega^*\cap\Omega)$ such that 
$u_n\rightarrow u_\infty$ in $L^2(\omega^\ast\cap\Omega)$. We show first that
$u_\infty$ is $A$-harmonic on $\omega^\ast\cap\Omega$ and belongs to $H^1_{loc}(\omega^\ast\cap\Omega)$. 
To see this pick any ball $B(x_0,r)\subset\subset\omega^*\cap\Omega$ centered at $x_0$ of radius $r$.
We apply the Cacciappoli inequality (Theorem \ref{Theorem 3.1}) together with the Rellich-Kondrachov compactness 
theorem to deduce that $u_n$ is Cauchy with respect to the energy norm in $B(x_0,r/2)$. From the completeness of $H^1(B(x_0,r/2)$ we see that $u_n\rightarrow u_\infty$ in $H^1(B(x_0,r/2))$. From this we conclude that $u_\infty\in H^1_{loc}(\omega^*\cap\Omega)$. The weak formulation
of the boundary value problem together with the strong convergence of the
sequence easily shows that $u_\infty$ is $A$-harmonic.

Next consider any open subset $\mathcal{O}\subset\omega^*\cap\Omega$ such
that $dist(\partial\mathcal{O},\partial\omega^*\cap\Omega)>0$ and $\partial%
\mathcal{O}\cap(\partial\Omega\cap\omega^*)\not=\emptyset$. Consider any
ball $B(x_0,r)$ centered at $x_0\in\partial\Omega$ of radius $r$ with $%
B(x_0,r)\cap\Omega$ contained inside $\mathcal{O}$. The dimorphism $\mathbf{x%
}=\mathbf{\Psi}^{-1}(\mathbf{y})$ maps $B(\mathbf{x}_0,r)\cap\Omega$ onto $%
\widetilde{B}(\mathbf{x}_0,r)\cap\widetilde{\Omega}$, with $\widetilde{B}(%
\mathbf{x}_0,r)\cap\partial\widetilde{\Omega}$ being part of the plane $%
y_n=0 $. Extend $\widetilde{B}(\mathbf{x}_0,r)\cap\widetilde{\Omega}$ across 
$y_n=0 $ by reflection. Denote the image of this extension under $\mathbf{%
\Psi}$ by $\mathcal{C}^{\ast}_E$ and define $\mathcal{C}^\ast=(B(x_0,r)\cap%
\Omega)\cup\mathcal{C}^{\ast}_E$. Now consider $u_n\in{H}_{A,0}^0(\omega^{\ast
}\cap\Omega)$ for which $u_n\rightarrow u_\infty$ in $L^2(\omega^{\ast
}\cap\Omega)$. Now $u_n$ can be extended as an $A-$harmonic function over $%
\mathcal{C}^\ast$ with $\Vert u_n\Vert_{L^2(\mathcal{C}^*)}\leq C\Vert
u_n\Vert_{L^2(B(x_0,r)\cap\Omega)}$ where $C$ depends only on $%
\partial\Omega $. Since $B(x_0,r/2)\cap\Omega\subset\mathcal{C}^\ast$ we
can apply a Cacciappoli inequality analogous to Theorem \ref{Theorem 3.1} to discover that
$\{u_n\}_{n=1}^\infty$ is a Cauchy sequence in $H^1(B(x_0,r/2)\cap\Omega)$
and we conclude that $u_n\rightarrow u_\infty$ in $H^1(B(x_0,r/2)\cap\Omega)$.
This establishes property (1) of the theorem.

Observe that since $\partial_\nu u_n\equiv n\cdot A\nabla u_n$ vanishes on $B(\mathbf{x}_0,r)\cap\partial\Omega$ we can write 
\begin{eqnarray}
&&\Vert \partial_\nu u_\infty \Vert_{H^{-1/2}(B(x_0,r/2)\cap\partial\Omega)}= \Vert
\partial_\nu u_n-\partial_\nu u_\infty\Vert_{H^{-1/2}(B(x_0,r/2)\cap\partial\Omega)}  \notag \\
&&\leq C
\Vert u_n-u_\infty\Vert_{H^{1}( B(x_0,r/2)\cap\Omega)}  \label{l2upper}
\end{eqnarray}
where $C$ is independent of $n$. Property (2) now follows on noting that 
\begin{eqnarray}
\lim_{n\rightarrow\infty}\Vert u_n-u_\infty\Vert_{H^{1}(
B(x_0,r/2)\cap\Omega)}=0.  \label{limitbdry}
\end{eqnarray}

\subsection{Nearly exponential convergence for GFEM applied to heterogeneous
systems}

\label{Macroexponential}

Theorems \ref{boundnearlyexpo} and \ref{boundnearlyexpoBDRY} provide the
local finite dimensional subspaces required for a global Galerkin
approximation with error that converges nearly exponentially with the
degrees of freedom. For a given partition $\omega_i$ , $i=1,\ldots,m$ we
denote these subspaces by $\hat{\Psi}_i$, $i=1,\ldots,m$. We augment each of the
subspaces $\hat{\Psi}_i$ with the subspace of constant functions denoted by $%
\mathbb{R}$ and write $\Psi_i=\hat{\Psi}_i\oplus\mathbb{R}$. 
For domains $\omega_i$ that touch the boundary of $\Omega$ the local approximations are taken from the hyperplane ${\Psi}_i\oplus u_p^i$. Here $u_p^i\in H_A(\omega^\ast_i\cap\Omega)$ is the local particular solution
introduced in the previous section that
satisfies $n\cdot A\nabla u_p=g$ on $\omega^\ast_i\cap\partial\Omega$. We denote the
dimensions of ${\Psi}_i$ by $N_i$ and set $N=\sum_{i=1}^m N_i$. 
Recalling Theorem \ref{globallocal} and applying Theorems \ref{boundnearlyexpo}
and \ref{boundnearlyexpoBDRY} we obtain the following approximation theorem.

\begin{theorem}
\label{Nearlyexpofem} Nearly exponential approximation for GFEM.\newline
For $\varepsilon>0$ there is a $N_\varepsilon>0$ such that for all $%
N>N_\varepsilon$, there exist $\zeta_i\in{\Psi}_i$, for $\omega_i\cap\partial\Omega=\emptyset$,
$\zeta_i\in{\Psi}_i\oplus u_p^i$, for $\omega_i\cap\partial\Omega\not=\emptyset$ and a
constant $\mathcal{K}$ independent of $N$ such that the approximation 
$\zeta\in H^1(\Omega)$ given by 
\begin{eqnarray}  \label{2.5opt}
\zeta(x)=\sum_{i=1}^{m}\zeta _{i}(x)\phi _{i}(x)
\end{eqnarray}
satisfies 
\begin{eqnarray}  \label{2.6bopt}
\parallel u_{0}-\zeta \parallel _{L^{2}(\Omega )} &\leq&\mathcal{K}
e^{-N^{\left(\frac{1}{1+d}-\epsilon\right)}}
\end{eqnarray}
and 
\begin{eqnarray}  \label{2.6copt}
\parallel u_{0}-\zeta \parallel _{\mathcal{E}(\Omega )}&\leq& \mathcal{K}
e^{-N^{\left(\frac{1}{1+d}-\epsilon\right)}}.
\end{eqnarray}
\end{theorem}

\setcounter{equation}{0} \setcounter{theorem}{0} \setcounter{lemma}{0}

\section{Implementation of the multiscale GFEM method}

\label{implementation}

In this section we provide an overview of the main ideas noting that the
specific challenges and details of the implementation are the focus of
future work. The implementation consists of three parts: 
\begin{itemize}
\item Multiple independent parallel computations for construction of the local bases and the subsequent assembly of the global stiffness matrix.
\item A single global computation using the global stiffness matrix and load vector.
\item
Recovery of preselected local features of the solution through the multiplication of the local bases by solutions of the global problem, e.g., the recovery of stresses at fiber matrix interfaces.
\end{itemize}

We now give an outline of the primary issues involved in the computation of the local 
optimal approximation spaces $\Psi _{i}$, $i=1,\cdots,m$ provided by
Theorems \ref{boundnearlyexpo} and \ref{boundnearlyexpoBDRY}.
The local bases are given by the eigenfunctions of the problems \eqref{concreteevalue} and %
\eqref{concreteevalueBDRY}. In what follows we will assume that all the subdomains $\omega_i$
are roughly the same size and we will suppress the index $i$ and write $\omega$ and $\Psi$.

A suitable and effective numerical method for the
construction of  the local basis $\Psi$ is given by the subspace approach (see \cite{106} chapter 11). This
method is based on the Raleigh-Ritz approximation. Here the key ingredient for the success of this method is 
the selection of a suitable subspace with span that should approximate the span of the
first $N_i$ eigenfunctions.  We now briefly discuss the construction of the subspace and the 
discrete representation of the eigenfunctions
used in the local basis $\Psi$ of dimension $N>1$.  We start by introducing $M>N$
functions $\varsigma _{k}$ , $k=1,2,..$ defined on the boundary $\partial \omega^{\ast }.$ 
The example given at the end of this section shows  that good
candidates for $\varsigma _{k}$ are the normal derivatives of the harmonic polynomials of
degree $k$. Other choices are also proposed in section 6. These
functions are then used to construct the $M$ dimensional subspace $\mathfrak{S}$\ of
A-harmonic functions $u_{\varsigma_k}$ on $\omega ^{\ast \text{ }}$ which
satisfy the Neumann boundary condition $n\cdot A\nabla u=\varsigma _{k}$.
This subspace $\mathfrak{S}$ is  used within the subspace approach to construct the desired $N$
eigenfunctions and eigenvalues. These
eigenfunctions will comprise the local approximation $\Psi$ used in the multiscale GFEM.
The appropriate selection of $N$ and $M$ is determined by the rate of decay
of the eigenvalues with respect to these parameters. The numerical construction
of the local basis for each subdomain can be carried out in
parallel.

Since we have established nearly exponential convergence it is expected that only a small number of
eigenfunctions will need to be computed. Moreover the examples presented in subsequent sections show that 
the desired eigenfunctions can be smooth on $\partial \omega
^{\ast }$ allowing for an accurate approximation for relatively small values of $M>N$.

The global basis is constructed by combining the local bases with the partition of unity functions. 
The partition of unity structure guarantees a sparse global stiffness matrix
and the assembly of the global stiffness matrix is also local procedure that can be carried
out in parallel.

We now provide a rough estimate for the computational work involved in the
GFEM for heterogeneous systems for problems posed over a computational
domain $\Omega \subset R^{2}$. To start we cover $\omega ^{\ast }$ by a finite element mesh with
elements of size $h$ and solve the Neumann problem for the differential
equation $div(A\nabla u)=0$, subject to the boundary condition $\varsigma _{k}$.
These solutions deliver the (approximate) A-harmonic functions $u_{\varsigma
}^{h}$. For example using Gaussian elimination (LU decomposition) we require
($\mid \omega ^{\ast }\mid h^{-2}$ )$^{-2}$ operations (due to the
sparseness of the stiffness matrix). Since $M<<h^{-2}$ the cost of the
computing $M$ functions does not change the order of operations. Furthermore
computation of the $N$ eigenfunctions is relatively small and hence the cost
of the creation of the space $\Psi$ is of order $\mid \omega _{i}\mid
^{-2}h^{-4}$. The computation of the entries in the associated stiffness
matrix will not change the order of operations. The major problem is the choice
of mesh size $h$ that leads to an acceptable accuracy. If the boundary
functions $\varsigma _{k}$ are smooth, for example the traces of harmonic
polynomials as mentioned above, then we conjecture that the
accuracy is on the order of $\varepsilon _{\varsigma_k}$=$\parallel
u_{\varsigma_k}-u_{\varsigma_k}^{h}\parallel _{\mathcal{E(\omega )}}\leq
Ch^{\gamma }M^{\delta }\parallel u\parallel _{\mathcal{E(\omega }^{\ast })}$. 
In future work we plan to analyze this conjecture. For the case of a proper mesh applied to the fiber
material discussed in section 2, we  expect $\gamma =2$ and the effect of 
$M$ is negligible. The span of the functions $u_{\varsigma_k}$ are then used to construct the approximate
space $\Psi ^{h}$ using the subspace method. These approximate spaces are then used within the GFEM
scheme.

On applying Gaussian elimination to the global stiffness matrix we
obtain the approximate solution over $\Omega$. The local error incurred by
approximating the solution $u$ over a subdomain $\omega $ using the
approximate local basis $\Psi^{h}$ will be of order $e^{-N^{\left(\frac{1}{1+d}-\epsilon\right)}}$. 
The size of the global stiffness matrix is 
$N\times\mid \Omega \mid /\mid \omega \mid$. Because the stiffness matrix is
sparse, the work of solving the global system using Gaussian elimination is $%
(N\times\mid\Omega \mid /\mid \omega \mid )^{-4}$. The final implementation issue involves the ``best choice'' of
the size of the local domains $\omega$ and $\omega^*$ for maximum  computational efficiency within the context of parallel computation. This question is not addressed here however it is clear that it is a very important issue that needs to be addressed in the implementation.

%Here we assume a weight characterizing the computational time when looking to ascertain the optimal size of the domains $%
%\omega$ as well the optimal width of the buffer $\omega ^{\ast }\backslash
%\omega^*$.

\begin{figure}[ht]
%[t]
\centering
\scalebox{0.3}{\includegraphics{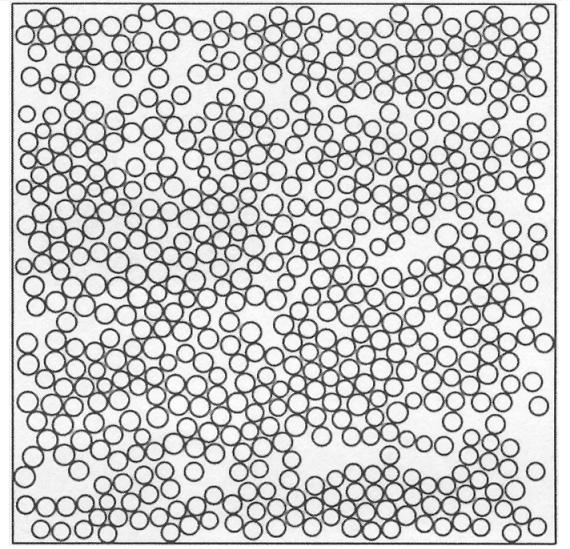}}
\caption{Mid-plane of shaft. Reprinted with permission from \cite{108}. Copyright 2004, John Wiley and Sons.}
\label{compdomain}
\end{figure}

We conclude with an example that illustrates the exponential convergence. 
In \cite{108} the numerical solution for the deformation inside a shaft reinforced with long  compliant fibers
with zero rigidity is given. The material in between the fibers is referred to as the matrix. Here the shaft is subjected to anti-pane shear loading and the system of elasticity reduces to the single scalar equation for deformations $u$ perpendicular to the mid-plane of the shaft.  When no fibers are in contact with each other  the entire
theory developed here also implies the exponential convergence of the GFEM for this problem. The computational domain is the shaft mid-plane $\Omega$ given by a subset of the $x-y$ plane portrayed in Figure \ref{compdomain}. For this problem when one constructs the local basis over a generic $\omega^*$ the convention is to remove any fiber domains intersecting $\partial\omega^*$ and to replace with matrix material.
For this kind of problem the local basis functions $u_{\varsigma _{k}}$ are taken to be harmonic in the matrix outside the fibers, taking zero Neumann data on the boundary of the fibers and taking Neumann data on the  boundary of $\partial\omega^*$ given by the traces of
harmonic polynomials of degree $k$. 
The Neumann condition  $\partial u/\partial n=g$ posed on the boundary of the computational domain $\partial\Omega$ is given by $g=(2x-y).$ We now describe the $\omega_i$ comprising the partition of unity for this example. The computational domain $\Omega$ is covered
by a $16\times 16$ mesh of square elements $\tau.$ The partition of unity
functions are the standard "hat" functions of the finite elements which are
bilinear on the elements $\tau.$ The supports of these "hat" functions create
the set of local domains $\omega_i.$  Sets $\omega_i $ interior to $\Omega$ are
composed of $2 \times 2$ squares. The interior domains $\omega_i ^{\ast }$ are
composed of  $6 \times 6$ squares. The domains $\omega_i $ and $\omega_i^{\ast }$
close to the boundary of $\Omega $ are constructed as described above in 
section 3.2.
The spaces $\Psi_i$ are the restrictions to $\omega_i\subset\omega_i^*$  of the solutions $u_{\varsigma _{k}}$ on $%
\omega_i ^{\ast }$ with  Neumann boundary condition given by the traces of
the harmonic polynomials of degree $k$, hence the dimension of $\Psi_i$ is $2k+1.$ These local basis functions
are computed using the finite element method on sufficiently fine
mesh, so that the error is  negligible.

\begin{figure}[t]
%[ht]
\centering
\scalebox{0.3}{\includegraphics{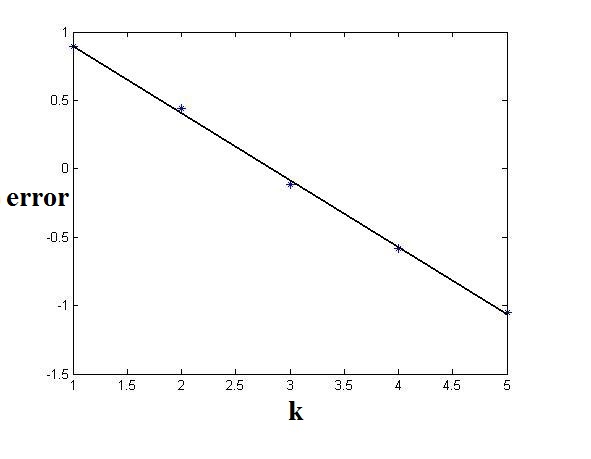}}
\caption{Decay of the error.}
\label{graph}
\end{figure}

The "exact" solution $u$ is computed by "overkill." The
relative energy norm of the error as function of the degree $k$ of the
harmonic polynomials is presented in Table 1. Figure \ref{graph} shows the error in the log
scale as function of $k$ in the linear scale. The straight line clearly
shows an exponential rate of decay. From these numbers we see that the rate is given by $e^{-0.48n}$
while the estimate is given by $e^{-0.33n},$  where $n=2k+1.$


The simulations presented in \cite{108} also numerically investigate the effect
of the distance between the boundaries of $\omega $ and $\omega ^{\ast }$.
There it is found that the rate of convergence decreases as the distance is
reduced and that the exponential convergence vanishes when $\omega$ and $%
\omega^*$ coincide.

\begin{equation*}
\begin{array}{|l|l|l|l|l|l|}
\hline
k & 1 & 2 & 3 & 4 & 5 \\ \hline
& $2.45\%$ & $1.55\%$ & $0.89\%$ & $0.56\%$ & $0.35\%$ \\ \hline
\end{array}%
\end{equation*}

\begin{center}
Table 1. The error as function of the degree of harmonic polynomials
charecterizing the boundary conditions $\varsigma _{\substack{ k  \\ }}$.
\end{center}

So far we addressed only the case that the right hand side $f=0$ in the
equation \eqref{1.1}. The general case $f\neq 0$ can be easily reduced to the case 
$f=0.$ To this end we introduce the local particular solution of the
differential equation \eqref{1.1} on each $\omega_i ^{\ast }$ subject to a constant  Neumann boundary condition over $\partial\omega_i^*$ determined according to the consistency condition. This is inexpensive to implement because the stiffness matrix has already been constructed together with its LU decomposition.
We denote the local particular solution by $u_{\omega_{\scriptscriptstyle{i}}}$. Then the local approximation over $\omega_i$ used in the GFEM belongs to the hyperplane  $\hat{\Psi}_i\oplus u_{\omega_{\scriptscriptstyle{i}}}$. Here the finite dimensional subspace $\hat{\Psi}_i $ is  given by the optimal local bases constructed for the $A-$harmonic problem.   As before this construction delivers a nearly exponential rate of convergence.

\setcounter{equation}{0} \setcounter{theorem}{0} \setcounter{lemma}{0}

\section{Homogenization of the $n$-width and exponential decay of
approximation error in the pre-asymptotic regime}

\label{nwidthhomog} 
%homogenization theorem for the eigenvalue problem \eqref{concreteevalue}.
We identify the homogenization limit of the $n$-width and the corresponding
optimal basis functions. These ideas are used to provide examples of
exponential convergence of the approximation error when the characteristic
length scale describing a heterogeneous medium is sufficiently small. In
what follows we work in the general context and homogenization is described
by $H$-convergence \cite{MuratTartar} or $G$-convergence \cite{7}. We
consider a sequence of coefficient matrices $A^\epsilon(x)$ in $\mathfrak{C}$
indexed by $\epsilon$, with $\epsilon=1/\ell$ for $\ell=1,2,\ldots$. Since
we consider symmetric coefficient matrices the notions of $G$ convergence
and $H$ convergence coincide and the class of coefficients $\mathfrak{C}$ is
compact with respect to $H$-convergence see \cite{MuratTartar}, \cite{7}. In
what follows we assume that the sequence $A^\epsilon$ $H$-converges to a
homogenized coefficient matrix $A^0$ in $\mathfrak{C}$ and we write $%
A^\epsilon\overset{H}{\rightarrow}A^0$.

We describe the $n$-widths associated with the sequence $A^\epsilon$ and the 
$H$-limit $A^0$. For each $A^\epsilon$ we introduce the Hilbert space $%
\mathcal{H}_\epsilon$ defined to be all elements in $H_{A^\epsilon}(%
\omega^*)/\mathbb{R}$ equipped with the energy inner product 
\begin{eqnarray}
(u,v)_{\mathcal{H}_\epsilon}=\int_{\omega^* }A^\epsilon\nabla u\cdot\nabla v
\, dx  \label{energyeps}
\end{eqnarray}
and norm $\Vert v\Vert_{\mathcal{H}_\epsilon}^2=(v,v)_{\mathcal{H}_\epsilon}$%
. The Hilbert space associated with $A^0$ is denoted by $\mathcal{H}_0$ and
is defined to be all elements in $H_{A^0}(\omega^*)/\mathbb{R}$ equipped
with the energy inner product 
\begin{eqnarray}
(u,v)_{\mathcal{H}_0}=\int_{\omega^* }A^0\nabla u\cdot\nabla v \, dx
\label{energyzed}
\end{eqnarray}
and norm $\Vert v\Vert_{\mathcal{H}_0}^2=(v,v)_{\mathcal{H}_0}$.

For each $A^\epsilon$ we introduce the restriction operator $P_\epsilon:%
\mathcal{H}_\epsilon\rightarrow H_{A^\epsilon}(\omega)/\mathbb{R}$ such that $%
P_\epsilon (u)(x)=u(x)$ for all $x\in\omega$ and $u\in \mathcal{H}_\epsilon$%
. As mentioned in the previous section the operator $P_\epsilon$ associated
with $A^\epsilon$ is compact. For future reference the energy bilinear form
defined on $H_{A^\epsilon}(\omega)/\mathbb{R}$ is given by 
\begin{eqnarray}
(u,v)_{\mathcal{E}_\epsilon(\omega)}=\int_{\omega}A^\epsilon \nabla
u\cdot\nabla v \, dx  \label{energyepsomeg}
\end{eqnarray}
and we set $\Vert u\Vert_{\mathcal{E}_\epsilon(\omega)}^2=(u,u)_{\mathcal{E}%
_\epsilon(\omega)}$. Similarly, for the $H$-limit $A^0$ we introduce the
compact operator given by the restriction $P_0:\mathcal{H}_0\rightarrow
H_{A^0}(\omega)/\mathbb{R}$ such that $P_0 (u)(x)=u(x)$ for all $x\in\omega$ and $u\in 
\mathcal{H}_0$. The energy bilinear form defined on $H_{A^0}(\omega)$ is
given by 
\begin{eqnarray}
(u,v)_{\mathcal{E}_0(\omega)}=\int_{\omega}A^0 \nabla u\cdot\nabla v \, dx
\label{energyzedomeg}
\end{eqnarray}
and we set $\Vert u\Vert_{\mathcal{E}_0(\omega)}^2=(u,u)_{\mathcal{E}%
_0(\omega)}$.

The $n$ width associated with the coefficients $A^\epsilon$ is given by 
\begin{eqnarray}
d_n^\epsilon(\omega,\omega^*)=\inf_{S(n)\subset
H_{A^\epsilon}(\omega)/\mathbb{R}}\sup_{u\in \mathcal{H}_\epsilon}\inf_{\chi\in S(n)}%
\frac{\Vert P_\epsilon u-\chi\Vert_{\mathcal{E_{\epsilon}}(\omega)}}{\Vert
u\Vert_{\mathcal{H}_\epsilon}}.  \label{nwidtheps}
\end{eqnarray}
The optimal local approximation space $\Psi^\epsilon_n(\omega)$ associated
with $A^\epsilon$ is described in terms of the eigenfunctions associated
with the following spectral problem. We introduce the adjoint operator $%
P_\epsilon^*:H_{A^\epsilon}(\omega)/\mathbb{R}\rightarrow \mathcal{H}_\epsilon$ and the
operator $P_\epsilon^*P_\epsilon$ is a self adjoint non-negative compact map
taking $\mathcal{H}_\epsilon$ into itself. The eigenfunctions and
eigenvalues are denoted by $\{\varphi^\epsilon_i\}_{i=1}^\infty$ and $%
\{\lambda^\epsilon_i\}_{i=1}^\infty$ and satisfy the problem 
\begin{eqnarray}
(P_\epsilon^*P_\epsilon\varphi_i^\epsilon ,\delta)_{\scriptscriptstyle{%
\mathcal{H}_\epsilon}}&=&\lambda^\epsilon_i (\varphi_i^\epsilon,\delta)_{%
\mathcal{H}_\epsilon}\hbox{,      $\forall \delta\in \mathcal{H}_\epsilon$}.
\label{vareigeneps}
\end{eqnarray}
The non-zero eigenvalues of $P_\epsilon^*P_\epsilon$ are listed according to
decreasing order of magnitude 
\begin{eqnarray}
\lambda_1^\epsilon\geq\lambda_2^\epsilon\geq\ldots>0.  \label{epsineq}
\end{eqnarray}
The optimal approximation space is given by $\Psi^\epsilon_n(\omega)=span\{%
\psi^\epsilon_1,\ldots,\psi^\epsilon_n\}$, where $\psi^\epsilon_i=P_\epsilon%
\varphi^\epsilon_i$ and $d_n^\epsilon(\omega,\omega^*)=\sqrt{%
\lambda_{n+1}^\epsilon}$. The $n$ width associated with the coefficient $A^0$
is given by 
\begin{eqnarray}
d_n^0(\omega,\omega^*)=\inf_{S(n)\subset H_{A^0}(\omega)}\sup_{u\in \mathcal{%
H}_0}\inf_{\chi\in S(n)}\frac{\Vert P_0 u-\chi\Vert_{\mathcal{E}_0(\omega)}}{%
\Vert u\Vert_{\mathcal{H}_0}}.  \label{nwidthzed}
\end{eqnarray}
The optimal local approximation space $\Psi^0_n(\omega)$ associated with $%
A^0 $ is described in terms of the eigenfunctions associated with the
following spectral problem. We introduce the adjoint operator $%
P_0^*:H_{A^\epsilon}(\omega)/\mathbb{R}\rightarrow \mathcal{H}_0$ and the operator $%
P_0^*P_0$ is a self adjoint non-negative compact map taking $\mathcal{H}_0$
into itself. The eigenfunctions and eigenvalues are denoted by $%
\{\varphi^0_i\}_{i=1}^\infty$ and $\{\lambda^0_i\}_{i=1}^\infty$ and satisfy
the problem 
\begin{eqnarray}
(P_0^*P_0 \varphi^0_i,\delta)_{\scriptscriptstyle{\mathcal{H}_0}%
}&=&\lambda^0_i (\varphi_i^0,\delta)_{\mathcal{H}_0}%
\hbox{,      $\forall
\delta\in \mathcal{H}_0$}.  \label{vareigenzed}
\end{eqnarray}
The non-zero eigenvalues of $P_0^*P_0$ are listed according to decreasing
order of magnitude 
\begin{eqnarray}
\lambda_1^0\geq\lambda_2^0\geq\ldots>0.  \label{zedineq}
\end{eqnarray}
The optimal approximation space is given by $\Psi^0_n(\omega)=span\{%
\psi^0_1,\ldots,\psi^0_n\}$, where $\psi^0_i=P_0\varphi^0_i$ and $%
d_n^0(\omega,\omega^*)=\sqrt{\lambda_{n+1}^0}$. The homogenization limit of $%
n$-widths and optimal approximations is given by the following theorem.

\begin{theorem}
Suppose that the coefficient matrices $A^\epsilon(x)$ in $\mathfrak{C}$ $H$%
-converge to $A^0(x)$ in $\mathfrak{C}$ as $\epsilon\rightarrow 0$. Then
there exists a subsequence of coefficients $A^\epsilon\overset{H}{\rightarrow%
}A^0$ such that 
\begin{eqnarray}
\lambda_i^\epsilon\rightarrow\lambda_i^0 &\hbox{ and }&
\varphi_i^\epsilon\rightharpoonup\varphi_i^0,%
\hbox{ weakly in
$H^1(\omega^*)$, for $i=1,2\ldots$}  \label{homognwidthevf}
\end{eqnarray}
Hence 
\begin{eqnarray}
\lim_{\epsilon\rightarrow
0}d_n^\epsilon(\omega,\omega^*)=d_n^0(\omega,\omega^*)
\label{homognwidnwidth}
\end{eqnarray}
and each function in the optimal basis for $A^\epsilon$ given by $%
\Psi^\epsilon_n(\omega)=span\{\psi^\epsilon_1,\ldots,\psi^\epsilon_n\}$
converges weakly in $H^1(\omega)$ to the corresponding function in the
optimal basis for $A^0$ given by $\Psi^0_n(\omega)=span\{\psi^0_1,\ldots,%
\psi^0_n\}$. \label{homognwidth}
\end{theorem}

%%%%%%%%%%%%%%%%%%%%%%%%%%%%%%%%
\noindent\emph{Proof.} The proof proceeds in three steps.

\noindent\textbf{Step 1}. We start by fixing the index $i$ and state the
following Lemma.

\begin{lemma}
Suppose $A^\epsilon\overset{H}{\rightarrow}A^0$. For $i$ fixed, consider the
associated eigenfunction and eigenvalue $(\varphi^\epsilon_i,\lambda_i^%
\epsilon)$ of $P_\epsilon^*P_\epsilon$. Sending $\epsilon\rightarrow 0$ and
passing to subsequences as necessary there exists a positive number $%
\overline{\lambda}_i$ and function $\overline{\varphi}_i\in\mathcal{H}_0$
for which 
\begin{eqnarray}
\lambda_i^\epsilon\rightarrow\overline{\lambda}_i &\hbox{ and }&
\varphi_i^\epsilon\rightharpoonup\overline{\varphi}_i, 
\hbox{weakly in
$H^1(\omega^*)$},  \label{ilim}
\end{eqnarray}
and 
\begin{eqnarray}
(P_0^*P_0 \overline{\varphi}_i,\delta)_{\scriptscriptstyle{\mathcal{H}_0}}&=&%
\overline{\lambda}_i (\overline{\varphi}_i,\delta)_{\mathcal{H}_0}%
\hbox{,   
$\forall \delta\in \mathcal{H}_0$}.  \label{vareigenoverline}
\end{eqnarray}
\label{homogevalueoverline}
\end{lemma}

Before proceeding with the proof of Lemma \ref{homogevalueoverline} we
state the following two compensated compactness results presented in the
work of Murat and Tartar \cite{MuratTartar} for later reference.

\begin{lemma}
\label{L1MuratTartar} Let $\epsilon\rightarrow 0$, $D$ be any open subset of 
$\mathbb{R}^d$, and $\xi^\epsilon$, $v^\epsilon$, be sequences such that 
\begin{eqnarray}
&&\xi^\epsilon\in L^2(D)^d,  \notag \\
&&\xi^\epsilon\rightharpoonup\xi^0\hbox{  weakly in $L^2(D)^d$},  \notag \\
&&div\xi^\epsilon\rightarrow div\xi^0\hbox{  strongly in $H^{-1}(D)$},
\label{Lemma1stuff}
\end{eqnarray}
and 
\begin{eqnarray}
&&v^\epsilon\in H^1(D),  \notag \\
&&v^\epsilon\rightharpoonup v^0\hbox{  weakly in $H^1(D)$}.
\label{Lemma1stuffbis}
\end{eqnarray}
Then 
\begin{eqnarray}
&&\int_D\xi^\epsilon\cdot\nabla v^\epsilon\eta\,dx\rightarrow
\int_D\xi^0\cdot\nabla v^0\eta\,dx, \hbox{  $\forall \eta\in C_0^\infty(D)$}.
\label{Lemma1stuffbis2}
\end{eqnarray}
\end{lemma}

\begin{lemma}
\label{T1MuratTartar} Suppose that $A^\epsilon(x)$ in $\mathfrak{C}$ $H$%
-converges to $A^0(x)$ in $\mathfrak{C}$ as $\epsilon\rightarrow 0$. Assume
that $u^\epsilon\in H^1(D)$, $f^\epsilon\in H^{-1}(D)$ and 
\begin{eqnarray}
-div(A^\epsilon\nabla u^\epsilon)&=&f^\epsilon,\hbox{  in $D$,}  \notag \\
u^\epsilon&\rightharpoonup& u^0 \hbox{  weakly in $H^1(D)$},  \notag \\
f^\epsilon&\rightarrow &f^0\hbox{  strongly in $H^{-1}(D)$},
\label{Th11stuff}
\end{eqnarray}
for $\epsilon\rightarrow 0$. Then 
\begin{eqnarray}
A^\epsilon\nabla u^\epsilon&\rightharpoonup& A^0\nabla u^0%
\hbox{  weakly in
$L^2(D)^d$,}  \notag \\
\int_D (A^\epsilon\nabla u^\epsilon\cdot\nabla
u^\epsilon)\eta\,dx&\rightarrow& \int_D (A^0\nabla u^0\cdot\nabla
u^0)\eta\,dx, \hbox{  $\forall \eta\in C_0^\infty(D)$  and}  \notag \\
A^\epsilon\nabla u^\epsilon\cdot\nabla u^\epsilon&\rightharpoonup& A^0\nabla
u^0\cdot\nabla u^0\hbox{  weakly in $L^1_{loc}(D)$}.  \label{Th11stuffbis2}
\end{eqnarray}
\end{lemma}

\noindent\emph{Proof of Lemma \ref{homogevalueoverline}.} Following
section 3 the we write \eqref{vareigeneps} as 
\begin{eqnarray}
\int_{\omega}A^\epsilon\nabla\varphi^\epsilon_i\cdot\nabla\delta\,dx=%
\lambda^\epsilon_i\int_{\omega^*}A^\epsilon\nabla\varphi^\epsilon_i\cdot%
\nabla\delta\,dx,\,\,\hbox{$\forall\delta\in \mathcal{H}_\epsilon$}.
\label{epsprob}
\end{eqnarray}

Now consider the sequence of eigenfunctions $\{\varphi^\epsilon_i\}_{%
\epsilon>0}$ for \eqref{epsprob} and with out loss of generality we
normalize $\varphi_i^\epsilon$ so that 
\begin{eqnarray}
\int_{\omega^*}
A^\epsilon\nabla\varphi_i^\epsilon\cdot\nabla\varphi_i^\epsilon\,dx=1&%
\hbox{
and   }&\lambda_i^\epsilon=\int_{\omega} A^\epsilon \nabla
\varphi_i^\epsilon\cdot\nabla\varphi_i^\epsilon\,dx.  \label{normal}
\end{eqnarray}
From \eqref{normal} we can extract a subsequence $\{\varphi_i^\epsilon\}_{%
\epsilon>0}$ and $\overline{\varphi}_i\in H^1(\omega^*)$ such that $%
\varphi_i^\epsilon\rightharpoonup\overline{\varphi}_i$ weakly in $%
H^1(\omega^*)$. Since $\varphi_i^\epsilon\in \mathcal{H}_\epsilon\subset
H^1(\omega^*)$ we apply Lemma \ref{T1MuratTartar} to deduce that for $%
\epsilon\rightarrow 0$, 
\begin{eqnarray}
A^\epsilon\nabla\varphi_i^\epsilon &\rightharpoonup & A^0\nabla\overline{%
\varphi}_i^,\hbox{ weakly in $L^2(\omega^*)^d,$}  \notag \\
\int_{\omega^*}(A^\epsilon\nabla\varphi_i^\epsilon\cdot\nabla\varphi_i^%
\epsilon)\delta\,dx &\rightharpoonup & \int_{\omega^*}(A^0\nabla\overline{%
\varphi}_i\cdot\nabla\overline{\varphi}_i)\delta\,dx,%
\hbox{$\forall
\delta\in C_0^\infty(\omega^*)$ and}  \notag \\
A^\epsilon\nabla\varphi_i^\epsilon\cdot\nabla\varphi_i^\epsilon
&\rightharpoonup & A^0\nabla\overline{\varphi}_i\cdot\nabla\overline{\varphi}%
_i,\hbox{ weakly in $L^1(\omega)$,}  \label{limits}
\end{eqnarray}
and $\overline{\varphi}_i \in \mathcal{H}_0$, hence 
\begin{eqnarray}
\lambda_i^\epsilon\rightarrow \int_{\omega}A^0 \nabla \overline{\varphi}%
_i\cdot\nabla\overline{\varphi}_i\,dx=\overline{\lambda}_i.
\label{limitlambda0}
\end{eqnarray}

To finish the proof we show that $\overline{\varphi}_i$ and $\overline{%
\lambda}_i$ are solutions of \eqref{vareigenoverline}. Consider any $g\in
H^{1/2}(\partial\omega^*)$ and $\delta^\epsilon\in \mathcal{H}_\epsilon$
such that $\delta^\epsilon=g$ on $\partial\omega^*$. We write $%
\delta^\epsilon= w^\epsilon+v$ where 
\begin{eqnarray}
v=g \hbox{ on }\partial\omega^* &\hbox{ and }& \Delta v=0 \hbox{ in }%
\omega^*,  \label{laplace}
\end{eqnarray}
and $w^\epsilon\in H_0^1(\omega^*)$, where 
\begin{eqnarray}
-\mathrm{div}(A^\epsilon\nabla w^\epsilon)=\mathrm{div}(A^\epsilon\nabla v),
\label{div}
\end{eqnarray}
in $\omega^*$.

%For any coefficient $A^\epsilon\in \mathfrak{C}$ one has that $\nabla
%w^\epsilon$ enjoys the higher integrability given by the following theorem
%of N. Meyers (Theorem 1 of \cite{Meyers}).

For any sequence of coefficients $A^\epsilon\in \mathfrak{C}$ Theorem 1 of \cite{Meyers}  shows that the sequence $\nabla
w^\epsilon$ enjoys the higher integrability given by the following Lemma.

\begin{lemma}
There is an interval $Q>p>2$ such that for $v\in W^{1,p}(\omega^*)$ then 
\begin{eqnarray}
\sup_{\varepsilon>0}\{\Vert \nabla w^\epsilon\Vert_{L^p(\omega^*)}\}<\infty.
\label{higherreg}
\end{eqnarray}
Here the interval is independent of $\epsilon$ and depends only on $\omega^*$%
, $\alpha$ and $\beta$. \label{lemmaa}
\end{lemma}

Now consider a dense subset $S$ of $H^{1/2}(\partial\omega^*)$ such that $%
g\in S$ implies that the solution $v$ of \eqref{laplace} belongs to $%
W^{1,p}(\omega^*)$. Then $\Vert\nabla\delta^\epsilon\Vert_{L^p(\omega^*)}%
\leq\Vert\nabla w^\epsilon\Vert_{L^p(\omega^*)}+\Vert\nabla
v\Vert_{L^p(\omega^*)}$ and Lemma \ref{lemmaa} implies that the associated
sequence $\{\delta^\epsilon\}_{\epsilon>0}$ satisfies 
\begin{eqnarray}
\sup_{\varepsilon>0}\{\Vert \nabla
\delta^\epsilon\Vert_{L^p(\omega^*)}\}<\infty.  \label{lpbound}
\end{eqnarray}
Additionally passing to a subsequence if necessary we see that there is an
element $\delta\in H^1(\omega^*)$ for which $\delta^\epsilon\rightharpoonup
\delta$ weakly in $H^1(\omega^*)$. Next an application of Lemma \ref%
{T1MuratTartar} shows that $\delta\in \mathcal{H}_0$ and an application of %
\ref{L1MuratTartar} to the sequences $\{A^\epsilon\nabla\varphi^\epsilon_i%
\}_{\epsilon>0}$ and $\{\delta^\epsilon\}_{\epsilon>0}$ gives 
\begin{eqnarray}
\int_{\omega^*}(A^\epsilon\nabla\varphi_i^\epsilon\cdot\nabla\delta^%
\epsilon)\eta\,dx &\rightarrow & \int_{\omega^*}(A^0\nabla\overline{\varphi}%
_i\cdot\nabla\delta)\eta\,dx \hbox{ $\forall \eta \in C_0^\infty(\omega^*)$}.
\label{morelimits}
\end{eqnarray}
From \eqref{lpbound} we deduce that $\{A^\epsilon\nabla\varphi_i^\epsilon%
\cdot\nabla\delta^\epsilon\}_{\epsilon>0}$ is equiintegrable on $\omega^*$
and it follows that 
\begin{eqnarray}
\int_{\omega}A^\epsilon\nabla\varphi_i^\epsilon\cdot\nabla\delta^\epsilon
&\rightarrow & \int_{\omega}A^0\nabla\overline{\varphi}_i\cdot\nabla\delta%
\,dx,  \label{morelimitsss1} \\
\int_{\omega^*}A^\epsilon\nabla\varphi_i^\epsilon\cdot\nabla\delta^\epsilon
&\rightarrow & \int_{\omega^*}A^0\nabla\overline{\varphi}_i\cdot\nabla\delta%
\,dx.  \label{morelimitsss2}
\end{eqnarray}
From \eqref{limitlambda0}, \eqref{morelimitsss1}, \eqref{morelimitsss2}, we
deduce that 
\begin{eqnarray}
\int_{\omega}A^0\nabla\overline{\varphi}_i\cdot\nabla\delta\,dx=\overline{%
\lambda}_i\int_{\omega^*}A^0\nabla\overline{\varphi}%
_i\cdot\nabla\delta\,dx,  \label{homogenizedprob}
\end{eqnarray}
or equivalently 
\begin{eqnarray}
(P_0^*P_0 \overline{\varphi}_i,\delta)_{\scriptscriptstyle{\mathcal{H}_0}}&=&%
\overline{\lambda}_i (\overline{\varphi}_i,\delta)_{\mathcal{H}_0}
\label{vareigenoverlinedense}
\end{eqnarray}
for all test functions $\delta$ belonging to $\mathcal{H}_0$ with traces in $%
S$, i.e., $\delta=g$ on $\partial\omega^*$, for $g\in S$. Lemma \ref%
{homogevalueoverline} now follows from the density of $S$ in $%
H^{1/2}(\partial\omega^*)$.

\noindent \textbf{Step 2}. We apply Lemma \ref{homogevalueoverline}
together with a diagonalization argument to extract a subsequence still
denoted by $(\lambda^\epsilon_i,\varphi^\epsilon_i)$, such that for every $%
i=1,2\ldots$ 
\begin{eqnarray}
\lambda^\epsilon_i\rightarrow\overline{\lambda}_i\hbox{ and }%
\varphi^\epsilon_i\rightharpoonup\overline{\varphi}_i%
\hbox{ weakly in
$H^1(\omega^*)$}  \label{diaglim}
\end{eqnarray}
where $(\overline{\lambda}_i,\overline{\varphi}_i)$ are solutions of %
\eqref{vareigenoverline} and 
\begin{eqnarray}
\overline{\lambda}_1\geq\overline{\lambda}_2\geq\ldots>0.
\label{listoverline}
\end{eqnarray}

\noindent\textbf{Step 3}. The final step is to show that all the
eigenfunctions and eigenvalues of the operator $P_0^*P_0$ are given by $(%
\overline{\lambda}_i,\overline{\varphi}_i)$ obtained in step 2. We argue by
contradiction and assume that there is an eigenvalue $\lambda$ of $P_0^*P_0$
for which $\lambda\not=\overline{\lambda}_i$ for every $i=1,2,\ldots$. Let $%
\varphi$ be a corresponding normalized eigenvector i.e., $%
(P_0^*P_0\varphi,\delta)_{\mathcal{H}_0}=\lambda(\varphi,\delta)_{\mathcal{H}_0}$,
for every $\delta\in\mathcal{H}_0$ and $\Vert\varphi\Vert_{\mathcal{H}_0}=1$%
. Then there is an integer $m$ such that 
\begin{eqnarray}
\lambda>\overline{\lambda}_{m+1}.  \label{hypcont}
\end{eqnarray}
To proceed we introduce the Rayleigh quotient for $v\in \mathcal{H}_\epsilon$
given by 
\begin{eqnarray}
R_\epsilon(v)=\frac{(P_\epsilon^*P_\epsilon v,v)_{\mathcal{H}_\epsilon}}{%
\Vert v\Vert^2_{\mathcal{H}_\epsilon}}  \label{rayleighquot}
\end{eqnarray}
and the eigenvalues of $P_\epsilon^*P_\epsilon$ listed in decreasing order
are given by 
\begin{eqnarray}
\lambda^\epsilon_i=\max_{v\in\mathcal{H}_\epsilon\perp\varphi_1^\epsilon,%
\ldots\varphi_{i-1}^\epsilon} R_\epsilon(v).  \label{max}
\end{eqnarray}
We establish the contradiction first under the extra assumption that the
gradient of $\varphi\in\mathcal{H}_0$ enjoys higher integrability and
belongs to $W^{1,p}(\omega^*)$ for $p>2$. We then indicate how to proceed
without this assumption.

Introduce $u^\epsilon\in\mathcal{H}_\epsilon$ such that $u^\epsilon=\varphi$
on $\partial\omega^*$. On passing to a further subsequence if needed we
apply Theorem \ref{T1MuratTartar} to see that 
\begin{eqnarray}
&&u^\epsilon\rightharpoonup\varphi\hbox{ weakly in $H^1(\omega^*)$ and} 
\notag \\
&&(P_\epsilon^*P_\epsilon u^\epsilon,u^\epsilon)_{\mathcal{H}%
_\epsilon}=\int_{\omega}A^\epsilon\nabla u^\epsilon\cdot\nabla
u^\epsilon\,dx\rightarrow\int_{\omega}A^0\nabla\varphi\cdot\nabla\varphi%
\,dx=\lambda.  \label{passagenum}
\end{eqnarray}
Noting from that $\varphi\in W^{1,p}(\omega^*)$ we observe from the
arguments preceding \eqref{lpbound} that 
\begin{eqnarray}
\sup_{\varepsilon>0}\{\Vert \nabla u^\epsilon\Vert_{L^p(\omega^*)}\}<\infty.
\label{lpbound2}
\end{eqnarray}
Thus the sequence $\{A^\epsilon\nabla u^\epsilon\cdot\nabla
u^\epsilon\}_{\epsilon>0}$ is equiintegrable on $\omega^*$ and we conclude
that 
\begin{eqnarray}
&&( u^\epsilon,u^\epsilon)_{\mathcal{H}_\epsilon}=\int_{\omega^*}A^\epsilon%
\nabla u^\epsilon\cdot\nabla
u^\epsilon\,dx\rightarrow\int_{\omega^*}A^0\nabla\varphi\cdot\nabla\varphi%
\,dx=1  \label{passagedenom}
\end{eqnarray}
so 
\begin{eqnarray}
\lim_{\epsilon\rightarrow 0}R_\epsilon(u^\epsilon)=\lambda.
\label{raylimit1}
\end{eqnarray}
Now introduce $v^\epsilon\in\mathcal{H}_\epsilon$ given by 
\begin{eqnarray}
v^\epsilon=u^\epsilon-\sum_{i=1}^m(u^\epsilon,\varphi^\epsilon_i)_{%
\scriptscriptstyle{\mathcal{H}_\epsilon}}\varphi^\epsilon_i  \label{proj}
\end{eqnarray}
As before we make use of the equiintegrability of $\{A^\epsilon\nabla
u^\epsilon\cdot\nabla\varphi_i^\epsilon\}_{\epsilon>0}$ on $\omega^*$
together with Lemma \ref{T1MuratTartar} to find that 
\begin{eqnarray}
(u^\epsilon,\varphi^\epsilon_i)_{\mathcal{H}_\epsilon}=\int_{\omega^*}A^%
\epsilon\nabla u^\epsilon\cdot\nabla
\varphi_i^\epsilon\,dx\rightarrow\int_{\omega^*}A^0\nabla\varphi\cdot\nabla%
\overline{\varphi}_i\,dx=(\varphi,\overline{\varphi}_i)_{\mathcal{H}_0}.
\label{componentlim}
\end{eqnarray}
Since $\lambda\not=\overline{\lambda}_i$ for all $i$, it follows from the
orthogonality of eigenvectors of $P_0^*P_0$ that $(\varphi,\overline{\varphi}%
_i)_{\mathcal{H}_0}=0$ for $i=1,2,\ldots,m$ and we deduce that 
\begin{eqnarray}
(u^\epsilon,\varphi^\epsilon_i)_{\mathcal{H}_\epsilon}\rightarrow 0.
\label{componentorthog}
\end{eqnarray}
Writing 
\begin{eqnarray}
&&\Vert v^\epsilon\Vert_{\mathcal{H}_\epsilon}^2=\Vert u^\epsilon\Vert_{%
\mathcal{H}_\epsilon}-\sum_{i=1}^m(u^\epsilon,\varphi^\epsilon_i)_{%
\scriptscriptstyle{\mathcal{H}_\epsilon}}^2,  \notag \\
&&(P_\epsilon^*P_\epsilon v^\epsilon,v^\epsilon)_{\mathcal{H}%
_\epsilon}=(P_\epsilon^*P_\epsilon u^\epsilon,u^\epsilon)_{\mathcal{H}%
_\epsilon}-\sum_{i=1}^m\lambda_i^\epsilon(u^\epsilon,\varphi^\epsilon_i)_{%
\scriptscriptstyle{\mathcal{H}_\epsilon}}^2  \label{sumss}
\end{eqnarray}
and sending $\epsilon$ to zero using \eqref{passagenum}, \eqref{passagedenom}%
, and \eqref{componentorthog} we conclude that 
\begin{eqnarray}
\lim_{\epsilon\rightarrow 0}R_\epsilon(v^\epsilon)=\lambda.  \label{contv}
\end{eqnarray}
On the other hand 
\begin{eqnarray}
(v^\epsilon,\varphi^\epsilon_i)_{\mathcal{H}_\epsilon}=0,%
\hbox{ for
$i=1,2,\ldots,m$}  \label{zeros}
\end{eqnarray}
so from \eqref{max} we get $\lambda^\epsilon_{m+1}\geq\lambda$ and taking
limits gives $\overline{\lambda}_{m+1}\geq\lambda$ which is a contradiction
to the original assumption $\lambda>\overline{\lambda}_{m+1}$.

We now remove the higher integrability assumption on the gradient of $%
\varphi\in\mathcal{H}_0$. For this case consider a sequence $%
s=1/\ell,\ell=1,2,\ldots$ and functions $\delta_s\in W^{1,p}(\omega^*)$  that converge to $\varphi$ in $%
W^{1,2}(\omega^*)$ as $s$ goes to zero. Choose $u_s^\epsilon\in\mathcal{H}%
_\epsilon$ such that $u^\epsilon_s=\delta_s$ on $\partial\omega^*$. Then
construct $v_s^\epsilon$ according to 
\begin{eqnarray}
v_s^\epsilon=u_s^\epsilon-\sum_{i=1}^m(u_s^\epsilon,\varphi_i^\epsilon)_{%
\mathcal{H}_\epsilon}\varphi^\epsilon_i.  \label{vscase}
\end{eqnarray}
As before $(v_s^\epsilon,\varphi_i^\epsilon)_{\mathcal{H}_\epsilon}=0$, for $%
i=1,\ldots,m$ and $\lambda_{m+1}^\epsilon\geq R_\epsilon(v_s^\epsilon)$.
Following previous arguments one deduces that the sequence $u_s^\epsilon$ is
bounded in $W^{1,p}(\omega^*)$ and on passing to subsequences as necessary 
\begin{eqnarray}
&&u_s^\epsilon\rightharpoonup u_s%
\hbox{ weakly in $H^1(\omega^*)$ where
$u_s\in\mathcal{H}_0$},  \notag \\
&&(P_\epsilon^*P_\epsilon u_s^\epsilon,u_s^\epsilon)_{\mathcal{H}%
_\epsilon}\rightarrow\int_{\omega}A^0\nabla u_s\cdot\nabla u_s\,dx,  \notag
\\
&&(u_s^\epsilon,u_s^\epsilon)_{\mathcal{H}_\epsilon}\rightarrow(u_s,u_s)_{%
\mathcal{H}_0},\hbox{ and}  \notag \\
&&(u_s^\epsilon,\varphi_i^\epsilon)_{\mathcal{H}_\epsilon}\rightarrow(u_s,%
\overline{\varphi}_i)_{\mathcal{H}_0},\hbox{ for $i=1,\ldots,m$}  \label{listepsilon}
\end{eqnarray}
and 
\begin{eqnarray}
\lim_{\varepsilon\rightarrow 0}R_\epsilon(v_s^\epsilon)=\frac{\int_\omega\,A^0%
\nabla u_s\cdot\nabla u_s dx-\sum_{i=1}^m\overline{\lambda}%
_i(u_s,\overline{\varphi}_i)_{\mathcal{H}_0}^2}{(u_s,u_s)_{\mathcal{H}%
_0}^2-\sum_{i=1}^m(u_s,\overline{\varphi}_i)_{\mathcal{H}_0}^2}.  \label{rs}
\end{eqnarray}
Since $\delta_s$ converges strongly in $H^1(\omega^*)$ to $\varphi$ it
follows from the uniqueness of solution of the Dirichlet boundary value
problem for $A^0$ harmonic functions that $u_s$ converges strongly to $%
\varphi$ in $\mathcal{H}_0$ thus 
\begin{eqnarray}
\overline{\lambda}_{m+1}\geq\lim_{s\rightarrow
0}\lim_{\varepsilon\rightarrow 0}R_\epsilon(v_s^\epsilon)=\lambda
\label{scontradiction}
\end{eqnarray}
and we arrive at a contradiction and Theorem \ref{homognwidth} is
proved.

We conclude by applying the homogenization of $n$-width theorem to construct
an example that shows exponential decay of the approximation error in the
pre-asymptotic regime. We consider a heterogeneous medium with
characteristic length scale $\epsilon>0$. To fix ideas we work in two
dimensions and suppose that the associated sequence of coefficients $%
A^\epsilon$ is such that it $H$-converges to a constant effective
conductivity $A^0$ matrix as $\epsilon\rightarrow 0$. In the coordinate
system corresponding to the eigenvectors $e^1$, $e^2$ of $A^0$ we have $%
A^0=a_1 e^1\otimes e^1+a_2e^2\otimes e^2$ and we set $b=a_2/a_1$. To fix
ideas we suppose that $\omega^*$ is the ellipsoid $E_{r^*}=%
\{(x_1,x_2);x_1^2+x_2^2/b=r^*\}$ and $\omega\subset\omega^*$ is the
concentric ellipsoid $E_r=\{(x_1,x_2);x_1^2+x_2^2/b=r\}$ with $r<r^*$. For 
$z=x+iy$ recall the harmonic polynomials $w_j(x_1,x_2)=\Re{z^j}=r^j\cos{%
(j\theta)}$, $\hat{w}_j(x_1,x_2)=\Im{z^j}=r^j\sin{(j\theta)}$, for $%
j=1,\ldots,n$. Calculation shows that the optimal basis associated with the $%
n$ width for $A^0$ is given by the $A^0$ harmonic polynomials $%
v_j=w_j(x_1,x_2/\sqrt{b})$, $\hat{v}_j=\hat{w}_j(x_1,x_2/\sqrt{b})$ and
eigenvalues $\lambda_j=e^{-2|\ln{\frac{r}{r^*}}|j}$, $j=1,\ldots,n$ of 
\begin{eqnarray}
\int_{E_r}A^0\nabla\varphi_j\cdot\nabla\delta\,dx=\lambda_j\int_{E_{r^*}}A^0%
\nabla\varphi_j\cdot\nabla\delta\,dx,  \label{eigenf}
\end{eqnarray}
for all $\delta\in H_{A^0}(\omega^*)$. It follows from Theorem \ref%
{thlambdanplusone} that the decay of approximation error for the optimal
basis associated with the homogenized coefficient $A^0$ is 
\begin{eqnarray}
e^{-|\ln{\frac{r}{r^*}}|(n+1)}.  \label{exponential}
\end{eqnarray}
Now we denote the $n$ width associated with the optimal basis for $%
A^\epsilon $ by $d_n^\epsilon(E_r,E_{r^*})$. Direct application of Theorem %
\ref{homognwidth} together with \eqref{exponential} gives the the following
bound on the pre-asymptotic rate of approximation error.

\begin{theorem}
Given $N>0$ and tolerance $\tau>0$ there exist an $\epsilon>0$ such that for 
$1\leq n\leq N$, that 
\begin{eqnarray}
e^{-|\ln{\frac{r}{r^*}}|(n+1)}-\tau\leq d_n^\epsilon(E_r,E_{r^*})\leq
e^{-|\ln{\frac{r}{r^*}}|(n+1)}+\tau.  \label{preestimate}
\end{eqnarray}
\label{preasymp}
\end{theorem}

\setcounter{equation}{0} \setcounter{theorem}{0} \setcounter{lemma}{0}

\section{Implementation in the pre-asymptotic regime and more examples of
exponential convergence}

\label{homogpre}

In this section we discuss a method for computational approximation that
employs the optimal basis for the homogenized problem to construct
approximation spaces for composites with heterogeneities on the length scale 
$\epsilon>0$ relative to the size of $\omega^*$. We work in the general
context and consider a sequence of coefficient matrices $\{A^\epsilon\}_{%
\epsilon>0}\in\mathfrak{C}$ that $H$-converge to a homogenized coefficient
matrix $A^0\in\mathfrak{C}$. For this case we recall the eigenfunctions $%
\varphi^\epsilon_i$ of \eqref{vareigeneps} and $\varphi^0_i$ of %
\eqref{vareigenzed} associated with $A^\epsilon$ and $A^0$ respectively. For 
$\epsilon>0$ fixed the optimal approximation space is given by the span of
the restriction of the functions $\varphi^\epsilon_i$, $i=1,\ldots,n$ to $%
\omega$. However in general it is known that the direct numerical
computation of eigenfunctions is computationally expensive. Instead we
introduce the functions $\phi_i^\epsilon\in H_{A^\epsilon}(\omega^*)/\mathbf{%
R}$ such that $\phi_i^\epsilon=\varphi_i^0$ on $\partial\omega^*$, for $%
i=1,\ldots,n$. We then define the approximation space $V_\epsilon^n(\omega)$
by 
\begin{eqnarray}
V_\epsilon^n(\omega)=span\left\{u^\epsilon_i=P\phi_i^\epsilon,\,\,\,i=1,%
\ldots,n \right\}  \label{apha}
\end{eqnarray}
and state the following approximation theorem

\begin{theorem}
\label{aphapprox} Given a tolerance $\tau>0$ there exists an $\overline{%
\varepsilon}>0$ such that $\forall\epsilon<\overline{\epsilon}$ 
\begin{eqnarray}
\Vert u_i^\epsilon-\varphi_i^\epsilon\Vert_{\mathcal{E}(\omega)}<\tau
\label{tolapprox}
\end{eqnarray}
\end{theorem}

We point out that this theorem remains the same if we choose $%
\phi_i^\epsilon\in H_{A^\epsilon}(\omega^*)/\mathbb{R}$ such that $n\cdot
A^\varepsilon\nabla\phi_i^\epsilon=n\cdot A^0\nabla\varphi_i^0$ on $%
\partial\omega^*$, for $i=1,\ldots,n$. When the homogenized coefficient $A^0$
is sufficiently simple e.g., $A^0$ is a constant, and $\omega$ and $\omega^*$
are concentric ellipsoids, the optimal approximation space for the
homogenized problem is given by explicit transcendental functions. And it
follows that the associated approximation space $V^n(\omega)$ is far less
expensive to compute than the eigenvalue problem associated with the optimal
approximation space. For these situations Theorem \ref{aphapprox} shows that 
$V_\epsilon^n(\omega)$ can be used provided that $\epsilon$ is sufficiently
small. We point out that the traces of the approximations $%
\phi_i^\varepsilon $ are indeed smooth on $\partial\omega^*$ noting that
this is exactly the assumption made in section when considering the accuracy
of the approximate local basis given in section \ref{implementation}. For
fiber reinforced composite materials it is clear that the size of $\omega$
needs to be chosen sufficiently large so that the relative length scale of
the fiber cross sections as characterized by $\varepsilon$ is sufficiently
small.

We now give the proof of Theorem \ref{aphapprox}. Recall from Theorem \ref%
{homognwidth} that $\varphi_i^\epsilon\rightharpoonup\varphi_i^0$ in $%
H^1(\omega^*)$, hence $\varphi_i^\epsilon\rightarrow\varphi_i^0$ in $%
L^2(\omega^*)$. On the other hand since $A^\epsilon$ $H$-converges to $A^0$
it follows from Theorem \ref{T1MuratTartar} that $\phi_i^\epsilon%
\rightharpoonup\varphi_i^0$ in $H^1(\omega^*)$, hence $\phi_i^\epsilon%
\rightarrow\varphi_i^0$ in $L^2(\omega^*)$. Application of the Caccioppoli
inequality delivers 
\begin{eqnarray}
&&\Vert u_i^\epsilon-\varphi_i^\epsilon\Vert_{\mathcal{E}(\omega)}<(4(\beta
)^{1/2}/\sigma\rho
)\Vert\phi_i^\epsilon-\varphi_i^\epsilon\Vert_{L^2(\omega^*)}  \notag \\
&&\leq(4(\beta )^{1/2}/\sigma\rho
)\left(\Vert\phi_i^\epsilon-\varphi_i^0\Vert_{L^2(\omega^*)}+\Vert\varphi_i^%
\epsilon-\varphi_i^0\Vert_{L^2(\omega^*)}\right)  \label{tolapproxa}
\end{eqnarray}
and Theorem \ref{aphapprox} is proved.

In the numerical example presented at the end of section \ref{implementation}
we have assumed that the homogenized equation is given by the Laplace
equation, i.e., $A^0=I$ and that the functions $\varphi_i^0=\varsigma_i$ are
the traces of harmonic polynomials.

Consider a family of heterogeneous media with characteristic length scale $%
\epsilon>0$. We suppose as before $A^\epsilon$ is $H$-convergent and
converges to a constant effective conductivity $A^0$ matrix as $%
\epsilon\rightarrow 0$. We take $\omega^*$ to be the unit square and $\omega$
to be a concentric square of side length $\sigma<1$ contained inside $%
\omega^*$. We suppose that $\sigma$ is such that we can fit concentric
ellipsoids $E_r\subset E_{r^*}$, with $r<r^*$ inside $\omega^*$ such that $%
\omega$ is contained inside the smaller ellipsoid $E_r$. We consider even
dimensional approximation spaces and take our approximation space $%
V_\epsilon^n(\omega^*)$ to be given by the span of the $A^\epsilon$-harmonic
functions $\phi_j$ on $\omega^*$ taking the Neumann data given by $%
\underline{n}\cdot A^0\nabla v_j$, for $j=1,\ldots n/2$ and $\hat{\phi}_j$
taking the Neumann data given by $\underline{n}\cdot A^0\nabla \hat{v}_j$,
for $j=1,\ldots, n/2$. Here $v_j$ and $\hat{v}_j$ are the $A^0$ harmonic
polynomials introduced in the previous section and $\underline{n}$ is the
outward directed unit normal on the boundary of $\omega^*$. For this case we
have the following theorem.

\begin{theorem}
For any sequence $\{u_\epsilon\}_{\epsilon>0}\subset
H_{A^\epsilon}(\omega^*)/\mathbb{R}$ such that $\sup_{\epsilon>0}\{\Vert
u_\epsilon\Vert_{\mathcal{E}_\epsilon(\omega^*)}\}<\infty$, then given $n>0$
and tolerance $\tau>0$ and on passing to a subsequence if necessary there
exist an $\epsilon_0>0$ such that for $\epsilon<\epsilon_0$ 
\begin{eqnarray}
\inf_{\chi\in V_{\epsilon}^n(\omega^*)}\Vert \chi-u_\epsilon\Vert_{\mathcal{%
E_\epsilon}(\omega)}\leq (e^{-|\ln{\frac{r}{r^*}}|(n+1)}+\tau)\Vert
u_\epsilon\Vert_{\mathcal{E_\epsilon}(\omega^*)}.  \label{preestimatesquare}
\end{eqnarray}
\label{preasympsquare}
\end{theorem}

\noindent\emph{Proof.} Let $\{\psi_1^\epsilon,\ldots,\psi_{n/2}^\epsilon,%
\hat{\psi}_1^\epsilon,\ldots,\hat{\psi}_{n/2}^\epsilon\}$ be the optimal
basis for the concentric ellipsoids $E_r\subset E_{r^*}$ for the coefficient 
$A^\epsilon$. The subspace spanned by these functions is denoted by $%
W_\epsilon^n(E_{r^*})$. The optimal basis for the concentric ellipsoids $%
E_r\subset E_{r^*}$ for the homogenized coefficient $A^0$, denoted by $%
W_0^n(E_{r^*})$, is precisely the span of the $A^0$ harmonic polynomials $%
v_j=w_j(x_1,x_2/\sqrt{b})$, $\hat{v}_j=\hat{w}_j(x_1,x_2/\sqrt{b})$, $%
j=0,\ldots,n/2$. Then there is a sequence of constant vectors $%
\{c_1^\epsilon,\ldots,c_{n/2}^\epsilon,\hat{c}_1^\epsilon,\ldots,\hat{c}%
_{n/2}^\epsilon\}$ bounded in $\mathbb{R}^n$ such that $\psi_\epsilon\in
W_\epsilon^n(E_{r^*})$ is given by $\psi_\epsilon=\sum_{j=1}^{n/2}
(c^\epsilon_j\psi_j^\epsilon+\hat{c}^\epsilon_j\hat{\psi}_j^\epsilon)$ and
for $\chi_\epsilon=\sum_{j=1}^{n/2} (c_j^\epsilon\phi^\epsilon_j+\hat{c}%
_j^\epsilon\hat{\phi}^\epsilon_j)$ we deduce that 
\begin{eqnarray}
&&\inf_{\chi\in V_{\epsilon}^n(\omega^*)}\Vert \chi-u_\epsilon\Vert_{%
\mathcal{E_\epsilon}(\omega)}\leq\Vert \chi_\epsilon-u_\epsilon\Vert_{%
\mathcal{E_\epsilon}(E_r)}  \notag \\
&&\leq\Vert u_\epsilon-\psi_\epsilon\Vert_{\mathcal{E_\epsilon}(E_r)} +\Vert
\chi_\epsilon-\psi_\epsilon\Vert_{\mathcal{E_\epsilon}(E_r)}  \notag \\
&&\leq d^n_\epsilon(E_r,E_{r^*})\Vert u_\epsilon\Vert_{\mathcal{E_\epsilon}%
(E_{r^*})}+\Vert \chi_\epsilon-\psi_\epsilon\Vert_{\mathcal{E_\epsilon}%
(E_{r})}  \notag \\
&&\leq d^n_\epsilon(E_r,E_{r^*})\Vert u_\epsilon\Vert_{\mathcal{E_\epsilon}%
(E_{r^*})}+ \frac{2(\beta )^{1/2}}{\delta} \Vert
\chi_\epsilon-\psi_\epsilon\Vert_{L^2(E_{r^*})}.  \label{limitsconcentric}
\end{eqnarray}
Here the last inequality in \eqref{limitsconcentric} follows from Theorem %
\ref{Theorem 3.1} and $\delta=dist(\partial E_{r^*},\partial E_r)$. Moreover
since $\{(c_1^\epsilon,\ldots,c_n^\epsilon,\hat{c}_1^\epsilon,\ldots,\hat{c}%
_{n/2}^\epsilon\}$ is bounded in $\mathbb{R}^n$ we can extract a convergent
subsequence and from our previous observations on $H$ convergence we have
that there is a $\psi_0\in W_0^n(E_{r^*})$ such that $\chi_\epsilon%
\rightarrow\psi_0$ and $\psi_\epsilon\rightarrow\psi_0$ strongly in $%
L^2(E_{r^*})$. It now follows that 
\begin{eqnarray}
&&\Vert \chi_\epsilon-u_\epsilon\Vert_{\mathcal{E_\epsilon}(\omega)}\leq
d^n_\epsilon(E_r,E_{r^*})\Vert u_\epsilon\Vert_{\mathcal{E_\epsilon}%
(E_{r^*})}  \notag \\
&&+\frac{2(\beta )^{1/2}}{\delta}\left(\Vert
\psi_\epsilon-\psi_0\Vert_{L^2(E_{r^*})}+\Vert
\chi_\epsilon-\psi_0\Vert_{L^2(E_{r^*})}\right),  \label{limitsconcentriclim}
\end{eqnarray}
and the theorem is proved.

Theorem \ref{preasympsquare} shows that the use of $V_{\epsilon}^n(\omega^*)$
delivers exponential convergence in the pre-asymptotic regime when the size
and separation of the disks is sufficiently small. 
\begin{figure}[ht]
%[t]
\centering
\scalebox{0.3}{\includegraphics{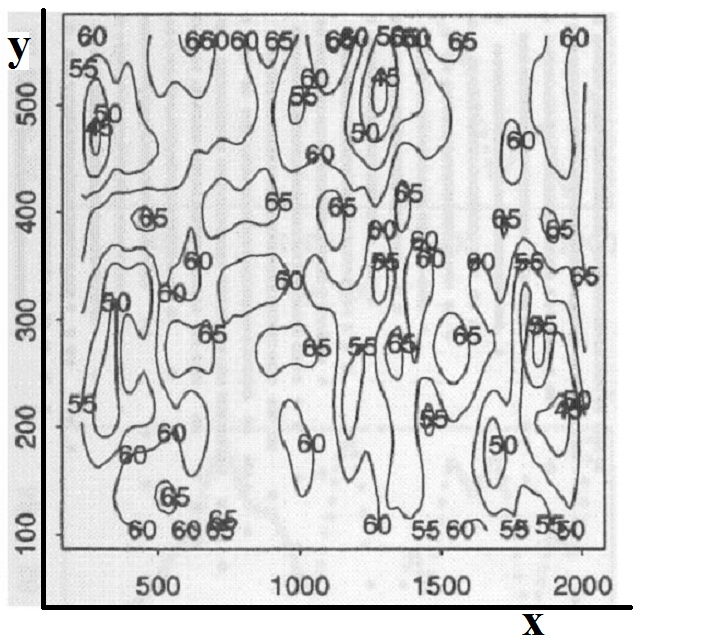}}
\caption{Level lines of fiber volume fraction.}
\label{level}
\end{figure}

These examples demonstrate how homogenized coefficients can be used in the
construction of the optimal shape functions. Of course the question of the
``best'' choice of homogenized coefficients which lead to a reasonable
approximation for general situations is not clear. Nevertheless one can
formally proceed by selecting appropriately sized $\omega ^{\ast }$ such
that it is large with respect to the features of the heterogeneity but such
that the heterogeneity is statistically uniform within it. With this in mind
we return to the fiber composite portrayed in Figure \ref{Plate}. Figure \ref%
{level} is a map of the level lines of the volume fraction taken over a
moving window given by the square of side length 116 $\mu m$, \cite{105}. It
provides a characterization of the spatial variability of the material. The
volume fraction varies between 45\% and 65\% across the sample and clearly
demonstrates the statistical inhomogeneity of the material. The correlation
between volume fraction and effective elastic properties for this sample is
illustrated in Figures 34 and 35 of \cite{105}. These Figures shows that the
spatial variation in effective properties correlates well with the variation
in volume fraction. With this in mind it appears that we should likely
choose $\omega^*$ to be of the size 200 $\mu m$.

We conclude this section by considering a periodic heterogeneous medium of
fixed period length given by $\epsilon=1/\ell>0$ where $\ell$ is a fixed
positive integer. Here we introduce a method for approximation that is
derived from the optimal basis associated with the homogenized coefficient
obtained from periodic homogenization. In what follows we will denote any
constant that is independent of $\epsilon$ and $n$ by $C$. We write the
coefficient describing the periodic medium $A^\epsilon(x)$ as a rescaling of
the coefficient of a unit periodic medium, i.e., $A^\epsilon(x)=A(x/%
\epsilon) $ where $A(y)$ is a coefficient of period one for $y\in \mathbb{R}%
^3$. We denote the unit period cell by $Q$ and the homogenized coefficient $%
A^0$ is given in terms of the periodic corrector matrix $P(y)$, $%
P_{ij}(y)=\partial_j w(y)^i+\delta_{ij}$ where $\chi(y)=(w^1,w^2,w^3)$ is
the $Q$ periodic solution of 
\begin{eqnarray}
div A(y)(\nabla \chi(y)+I)=0, \hbox{for y in Q},  \label{period}
\end{eqnarray}
and 
\begin{eqnarray}
A^0=\int_Q P(y)dy.  \label{effcoeff}
\end{eqnarray}

The optimal basis for $A^\epsilon$ is given in terms of the eigenfunctions $%
\varphi^\epsilon_i$ of \eqref{vareigeneps}. The optimal basis for the
homogenized coefficient with $A^0$ is given in terms of the eigenfunctions $%
\varphi^0_i$ of \eqref{vareigenzed}. We fix $\epsilon=1/\ell>0$ and the
optimal approximation space is given by the span of the restriction of the
functions $\varphi^\epsilon_i$, $i=1,\ldots,n$ to $\omega$. In this
implementation we introduce the functions $\phi_i^\epsilon\in
H_{A^\epsilon}(\omega^*)/\mathbb{R}$ such that $\phi_i^\epsilon=\varphi_i^0$
on $\partial\omega^*$, for $i=1,\ldots,n$. As before we define the
approximation space $V_\epsilon^n(\omega)$ by 
\begin{eqnarray}
V_\epsilon^n(\omega)=span\left\{u^\epsilon_i=P\phi_i^\epsilon,\,\,\,i=1,%
\ldots,n \right\}  \label{aphaH}
\end{eqnarray}
and state the following approximation theorem.

\begin{theorem}
\label{aphapproxH} Given any function $u\in H_{A^\epsilon}(\omega^*)/\mathbf{%
R}$ then 
\begin{eqnarray}
\min_{w\in V_\epsilon^n(\omega)}\left\{\Vert u-w\Vert_{\mathcal{E_\epsilon}%
(\omega)}\right\}\leq( d_{n-1}^\epsilon+ C\,\epsilon)\Vert u\Vert_{\mathcal{%
E_\epsilon}(\omega^*)},  \label{tolapproxH}
\end{eqnarray}
where $d_{n-1}^\epsilon$ is the $n$-width associated with $A^\epsilon$.
Moreover $d_{n-1}^\epsilon$ is estimated in terms of an easily computable
quantity 
\begin{eqnarray}
Q_\epsilon^n=\sqrt{\int_\omega A^\epsilon\nabla u^\epsilon_n\cdot\nabla
u^\epsilon_n\,dx}  \label{Q}
\end{eqnarray}
and the estimate is given by 
\begin{eqnarray}
|Q_\epsilon^n-d_{n-1}^\epsilon|\leq\,C\epsilon^{1/2}.  \label{dest}
\end{eqnarray}
\end{theorem}

\noindent\emph{Proof.} The theorem is proved by constructing upper bounds on
the quantity 
\begin{eqnarray}
R_\ell=\sup_{\{u\in H_{A^\epsilon};\,\Vert u\Vert_{\mathcal{E}%
_\epsilon(\omega^*)}=1\}}\left\{\inf_{w\in V^n(\omega)}\Vert u-w\Vert_{%
\mathcal{E}_\epsilon(\omega)}\right\}.  \label{boundupp}
\end{eqnarray}
From the corrector theory of periodic homogenization it follows from \cite%
{JKO} that there exists a constant $C$ depending only on $\Vert
D^2\varphi_i^0\Vert_{L^2(\omega^*)}$ and $\alpha<\beta$ for which 
\begin{eqnarray}
\Vert \varphi^\epsilon_i-\varphi^0_i\Vert_{L^2(\omega^*)}\leq \,C\epsilon,
\label{l2}
\end{eqnarray}
and since $\phi_i^\epsilon\in H_{A^\epsilon}(\omega^*)/\mathbb{R}$ with $%
A^\epsilon$ G-converging to $A^H$ and $\phi_i^\epsilon\rightharpoonup%
\varphi^0_i$ it follows again from \cite{JKO} that 
\begin{eqnarray}
\Vert \phi^\epsilon_i-\varphi^0_i\Vert_{L^2(\omega^*)}\leq \,C \epsilon.
\label{l2next}
\end{eqnarray}
Hence 
\begin{eqnarray}
\Vert \varphi^\epsilon_i-\phi^\epsilon_i\Vert_{L^2(\omega^*)}\leq
\,C\epsilon.  \label{tri}
\end{eqnarray}

Now consider $u\in H_{A^\epsilon}(\omega^*)/\mathbb{R}$ with $\Vert u \Vert_{%
\mathcal{E}_\epsilon(\omega^*)}=1$. For $w\in V_\epsilon^n(\omega)$ there
are constants $c_1,c_2,\ldots,c_n$ such that we can write $w=\sum_{i=1}^n
c_i u_i^\epsilon$ and we choose these constants $c_1,c_2,\ldots,c_n$ such
that $\varphi^\epsilon=\sum_{i=1}^n c_i \varphi_i^\epsilon$ gives the
optimal approximation to $u$ in the $\mathcal{E}_\epsilon(\omega)$ norm. For
this choice one has 
\begin{eqnarray}
\Vert u-w\Vert_{\mathcal{E}_\epsilon(\omega)}&\leq &\Vert u-\sum_{i=1}^n c_i
\varphi_i^\epsilon\Vert_{\mathcal{E}_\epsilon(\omega)}  \notag \\
&+&\Vert \sum_{i=1}^n c_i (\varphi_i^\epsilon-\phi^\epsilon)\Vert_{\mathcal{E%
}_\epsilon(\omega)}  \notag \\
&\leq&d_{n-1}^\epsilon+\,C\epsilon,  \label{startingineq}
\end{eqnarray}
where the first term on the last line of the inequality follows from the
definition of $n$-width and optimal basis and the second term follows from %
\eqref{tri} and it follows that $R_\ell\leq d_{n-1}^\epsilon+ C\epsilon$.

We conclude the proof by establishing \eqref{dest}. From Theorem \ref%
{thlambdanplusone} 
\begin{eqnarray}
(d_{n-1}^\epsilon)^2=\int_\omega
A^\epsilon\nabla\varphi_n^\epsilon\cdot\nabla\varphi_n^\epsilon\,dx,
\label{destineq}
\end{eqnarray}
where we have taken the normalization 
\begin{eqnarray}
\int_{\omega^*}
A^\epsilon\nabla\varphi_n^\epsilon\cdot\nabla\varphi_n^\epsilon\,dx=1.
\label{destineqnorm}
\end{eqnarray}
On choosing $u_n^\epsilon\in V_\epsilon^n(\omega)$ we write 
\begin{eqnarray}
(Q_\epsilon^n)^2=\int_\omega A^\epsilon\nabla u_n^\epsilon\cdot\nabla
u_n^\epsilon\,dx  \label{destineqq}
\end{eqnarray}
and 
\begin{eqnarray}
(Q_\epsilon^n)^2-(d_{n-1}^\epsilon)^2=\int_\omega A^\epsilon(\nabla
u_n^\epsilon+\nabla \varphi_n^\epsilon)\cdot(\nabla u_n^\epsilon-\nabla
\varphi_n^\epsilon)\,dx.  \label{destsumdiff}
\end{eqnarray}
Apriori elliptic estimates show that $C=\sup_{\epsilon>0}\{\Vert
u_n^\epsilon+\varphi_n^\epsilon\Vert_{\mathcal{E}_\epsilon(\omega^*)}\}<%
\infty$ and 
\begin{eqnarray}
|(Q_\epsilon^n)^2-(d_{n-1}^\epsilon)^2|&\leq& C \sqrt{\int_\omega
A^\epsilon(\nabla u_n^\epsilon-\nabla \varphi_n^\epsilon)\cdot(\nabla
u_n^\epsilon-\nabla \varphi_n^\epsilon)\,dx}  \notag \\
&\leq& C\,\Vert \phi_n^\epsilon- \varphi_n^\epsilon\Vert_{L^2(\omega^*)}
\leq \,C \epsilon  \label{destdiffbd}
\end{eqnarray}
where the second to last inequality follows from Theorem \ref{Theorem 3.1}
and the last inequality follows from \eqref{tri}. Inequality \eqref{dest}
follows noting that 
\begin{eqnarray}
|Q_\epsilon^n-d_{n-1}^\epsilon|\leq
|(Q_\epsilon^n)^2-(d_{n-1}^\epsilon)^2|^{1/2}.  \label{holder}
\end{eqnarray}

\setcounter{equation}{0} \setcounter{theorem}{0} \setcounter{lemma}{0} \appendix

\section{Appendix}

We provide a proof of the Cacciappoli inquality given in Lemma \ref%
{Theorem 3.1}. We introduce the cut off function $\eta\in C_0^1(\omega^\ast)$
such that $0\leq \eta\leq 1$ and $\eta=1$ for points inside $\mathcal{O}$
and $|\nabla\eta(x)|\leq 1/\delta$ for points in $\omega^\ast$. Given the
function $u\in \mathcal{H}(\omega^{\ast})$ and since $u$ is A -- harmonic we
have 
\begin{eqnarray}
\int_{\omega^\ast}\, A\nabla u\cdot\nabla(\eta^2 u)\,dx=0.  \label{identity}
\end{eqnarray}
Expanding (\ref{identity}) gives 
\begin{eqnarray}
&&\int_{\omega^\ast}\, (A\nabla u\cdot\nabla u) \eta^2\, dx =
-2\int_{\omega^\ast}\, (\eta A^{1/2}\nabla u)\cdot (u A^{1/2}\nabla\eta)\, dx
\notag \\
&&\leq 2\left(\int_{\omega^\ast}\, (A\nabla u\cdot\nabla u) \eta^2\,
dx\right)^{1/2}\left(\int_{\omega^\ast}\, (A\nabla \eta\cdot\nabla \eta)
u^2\, dx\right)^{1/2}  \label{cauchy}
\end{eqnarray}
so 
\begin{eqnarray}
&&\Vert u\Vert_{\mathcal{E}(\mathcal{O})}\leq \left(\int_{\omega^\ast}\,
(A\nabla u\cdot\nabla u) \eta^2\, dx\right)^{1/2} \leq
2\left(\int_{\omega^\ast}\, (A\nabla \eta\cdot\nabla \eta)u^2\,
dx\right)^{1/2}  \notag \\
&&\leq 2\gamma_2\left(\int_{\omega^\ast}\, |\nabla \eta|^2 u^2\,
dx\right)^{1/2} \leq\frac{2\gamma_2^{1/2}}{\delta}\Vert u
\Vert_{L^2(\omega^\ast)}.  \label{uppper}
\end{eqnarray}
and Lemma \ref{Theorem 3.1} is proved.

We now show that the restriction operators introduced in section three are
compact. We first consider two concentric cubes $\omega\subset\omega^\ast$.
The restriction operator $P:H_A(\omega^*)/\mathbb{R}\rightarrow H_A(\omega)/\mathbb{R}$ is defined
by $Pu(x)=u(x)$ for all $x\in\omega$ and all $u\in H_A(\omega^*)$.

\begin{lemma}
\label{compact} Given any sequence $\{u_n\}_{n=1}^\infty\in H_A(\omega^*)/\mathbb{R}$
that is bounded in the energy norm over $(\omega^*)$ then one can extract a subsequence that
converges in  $H^1(\omega)$ to an element of $H_A(\omega)/\mathbb{R}$.
\end{lemma}

\noindent\emph{Proof.} We apply the Poincare inequality together with the Rellich compactness theorem to extract a convergent
subsequence in $L^2(\omega^*)$. From Lemma \ref{Theorem 3.1} it now
follows that this subsequence is Cauchy with respect to the energy norm over 
$\omega$ and the convergence in $H^1(\omega)$ follows. The weak formulation
of the boundary value problem together with the strong convergence of the
subsequence easily shows that the limit function is $A$-harmonic and the
theorem is proved.

Next we consider two concentric cubes $C\subset\omega^*$ such that $%
\omega=C\cap\Omega$ and $\omega^*\cap\Omega$ have non zero volume. Here the
side length of $C$ is $\sigma$ and that of $\omega^*$ is $%
\sigma^*=(1+\rho)\sigma$. The restriction operator $P:H_{A,0}(\omega^*\cap%
\Omega)/\mathbb{R}\rightarrow H_{A,0}(\omega)/\mathbb{R}$ is defined by $Pu(x)=u(x)$ for all $%
x\in\omega$ and all $u\in H_{A,0}(\omega^*\cap\Omega)$. Here we suppose the
boundary of $\Omega$ is $C^1$.

\begin{lemma}
\label{CompactBDRY} Given any sequence $\{u_n\}_{n=1}^\infty\in
H_{A,0}(\omega^*\cap\Omega)/\mathbb{R}$ that is bounded with respect to the energy norm $(\omega^*\cap\Omega)$
then one can extract a subsequence that converges in  $H^1(\omega)$ 
to an element of $H_{A,0}(\omega)/\mathbb{R}$.
\end{lemma}

\noindent\emph{Proof.} Following section 3 we extend each $u_n\in
H^1_{A,0}(\omega^*\cap\Omega)/\mathbb{R}$ as an $A$-harmonic function across $\partial\Omega$
onto the set $\omega_E^*$  such that 
\begin{eqnarray}
\Vert u_n\Vert_{H^1(\omega_E^*)}\leq C \Vert
u_n\Vert_{H^1(\omega^*\cap\Omega)}  \label{extendcompact}
\end{eqnarray}
where $C$ depends only on $\partial\Omega$. Application of Theorem \ref%
{Theorem 3.1} gives 
\begin{eqnarray}
\Vert u_n\Vert_{\mathcal{E}(\omega)}\leq \frac{4\beta^{1/2}}{\sigma\rho}
\Vert u_n\Vert_{L^2(\omega_E^*)}  \label{cacciappoliappendix}
\end{eqnarray}
and we deduce that 
\begin{eqnarray}
\Vert u_n\Vert_{\mathcal{E}(\omega)}\leq C\frac{4\beta^{1/2}}{\sigma\rho}
\Vert u_n\Vert_{L^2(\omega^*\cap\Omega)}.  \label{cacciappoliappendixfinal}
\end{eqnarray}
With \eqref{cacciappoliappendixfinal} in hand we can now proceed as in the
proof of Lemma \ref{compact} to establish compactness.

\end{document}